\documentclass[11pt]{amsart}
\usepackage[utf8]{inputenc}
\usepackage[T1]{fontenc}
\usepackage{amsmath,amssymb,amsthm,mathrsfs}
\usepackage[a4paper,margin=3cm]{geometry}
\usepackage[colorlinks=true,linkcolor=blue,citecolor=blue]{hyperref}
\usepackage{booktabs}

\usepackage{xcolor}

\theoremstyle{plain}
\newtheorem{theorem}{Theorem}[section]
\newtheorem{proposition}[theorem]{Proposition}
\newtheorem{lemma}[theorem]{Lemma}
\newtheorem{corollary}[theorem]{Corollary}
\theoremstyle{definition}
\newtheorem{definition}[theorem]{Definition}
\newtheorem{example}[theorem]{Example}
\newtheorem{remark}[theorem]{Remark}

\newcommand{\g}{\mathfrak{g}}
\newcommand{\h}{\mathfrak{h}}
\newcommand{\q}{\mathfrak{q}}
\newcommand{\s}{\mathfrak{s}}
\newcommand{\aaa}{\mathfrak{a}}
\newcommand{\n}{\mathfrak{n}}
\newcommand{\dd}{\mathfrak{d}}
\newcommand{\kk}{\mathfrak{k}}
\newcommand{\m}{\mathfrak{m}}
\newcommand{\R}{\mathbb{R}}
\newcommand{\CH}{\mathbb{C}H}
\newcommand{\ip}[2]{\langle #1,#2\rangle}

\DeclareMathOperator{\ad}{ad}
\DeclareMathOperator{\Der}{Der}
\DeclareMathOperator{\Ker}{Ker}
\DeclareMathOperator{\Imm}{Im}

\DeclareMathOperator{\spann}{span}

\begin{document}

\title[$\varphi$-symmetric Sasakian Lie algebras]{$\varphi$-symmetric Sasakian Lie algebras}

\author{Adri\'an Andrada}
\address{(AA) FaMAF, Universidad Nacional de C\'ordoba and CIEM-CONICET,
Av.~Medina Allende S/N, Ciudad Universitaria, X5000HUA C\'ordoba (AR)}
\email{adrian.andrada@unc.edu.ar}

\author{Alberth J. N\'u\~nez Sullca}
\address{(AN) Facultad de Ciencias, Universidad Nacional de Ingenier\'ia,
Av.~T\'upac Amaru 210, R\'imac, Lima (PE)}
\email{alberth.nunez.s@uni.edu.pe}


\begin{abstract}
We study Sasakian Lie algebras with trivial centre for which the adjoint operator $\ad_\xi$ associated to the Reeb vector field $\xi$ has two-dimensional
image. Their structure is governed by a K\"ahler-exact Lie algebra $\h_1$ and a vector
$H_0\in\h_1$, and the cases $H_0=0$ and $H_0\neq0$ are of a different nature. For $H_0=0$ we
prove that the Sasakian structure is $\varphi$-symmetric if and only if $\h_1$ is locally
symmetric. For $H_0\neq0$ and $\h_1$ solvable, we prove
that it is $\varphi$-symmetric if and only if $\h_1$ is locally symmetric and a natural
two-form $\rho_0$ on $\h_1$ determined by $H_0$ is parallel. The latter condition is
equivalent to the vanishing of the mixed root spaces attached to the distinguished root
selected by $H_0$ in the associated normal $j$-algebra. As a consequence, every Sasakian Lie
algebra in this class of dimension at most seven is $\varphi$-symmetric, which extends a result of Loiudice and Lotta in dimension five. This bound is sharp: in dimension nine we exhibit a Sasakian structure which is not
$\varphi$-symmetric, even though its K\"ahler base is locally symmetric. 
\end{abstract}

\maketitle

\section{Introduction}

Sasakian manifolds are the odd-dimensional analogue of K\"ahler manifolds: a Riemannian
manifold $(M,g)$ of dimension $2n+1$ is Sasakian if and only if the cone $M\times\R^+$ is
K\"ahler. It is well known that local symmetry is  too restrictive in this context, since a locally symmetric Sasakian
manifold has constant curvature \cite{Ok}. Takahashi \cite{Ta} introduced instead the class
of \emph{$\varphi$-symmetric} spaces, in which the reflections along the integral curves of
the Reeb field are local automorphisms of the structure. Equivalently, the local K\"ahler
quotients by the Reeb foliation are Hermitian locally symmetric spaces. In \cite{BV1} the classification of 3-dimensional simply connected $\varphi$-symmetric spaces was given, while in \cite{BV2} it was proved that complete simply connected $\varphi$-symmetric spaces are necessarily naturally reductive homogeneous spaces. Later, the (local) classification in dimension five appeared in \cite{KW}, while the local classification in any dimension was given in \cite{JK}.

For left-invariant structures the question is which Lie groups carry a $\varphi$-symmetric Sasakian structure. In dimension five the answer is known. The case of nontrivial centre was settled by Calvaruso and Fino \cite[Theorem 5.3]{CF}, who showed that such a Lie algebra is $\varphi$-symmetric unless it is the contactization of the Kähler Lie algebra $\mathfrak d_{4,2}$, which is not locally symmetric. As for trivial centre, Loiudice and Lotta \cite[Theorem 5]{LL}  proved that every five-dimensional Sasakian Lie algebra with trivial centre is $\varphi$-symmetric. We do not know of any results beyond dimension five. This paper studies how far that statement extends.

Our setting is that of \cite{ACN}. If $\g$ is a Sasakian Lie algebra with trivial centre,
then $\ad_\xi$ is skew-symmetric and $\xi\in\Ker\ad_\xi$, so $\dim\Imm\ad_\xi$ is even and
positive; the first case is $\dim\Imm\ad_\xi=2$. By \cite[Proposition 4.2 and Theorem
4.3]{ACN}, such a $\g$ is built from a K\"ahler-exact Lie algebra $\h_1$ with Lie bracket denoted by $\theta(\cdot,\cdot)$, a vector
$H_0\in\theta(\h_1,\h_1)$ and data $(k,\gamma,\mu)$ subject to certain compatibility conditions.
Moreover, $\g$ is solvable if and only if $H_0\neq0$ and $\h_1$ is solvable
\cite[Proposition 4.5]{ACN}, so the case $H_0=0$ and the solvable case are different in
nature and must be treated with different approaches.

Write $\rho_0(x,y)=\ip{\theta(x,y)}{H_0}$ for the two-form on $\h_1$ that detects the
component of the bracket along $H_0$. We treat the two cases separately. For $H_0=0$ no
solvability assumption is needed and the answer is Theorem \ref{H0zero}, which states:
\[
\g \text{ is }\varphi\text{-symmetric}\iff \h_1 \text{ is locally symmetric,}
\]
that is, any Lie group with Lie algebra $\h_1$, equipped with the left-invariant metric induced by the inner product on $\h_1$, is locally symmetric. In this case $\g$ splits as a direct product of $\mathfrak{su}(2)$ or $\mathfrak{sl}(2,\R)$ with a
copy of $\h_1$, so the K\"ahler quotient of $G$ by the Reeb foliation is the Riemannian product of
a surface of constant curvature and the Lie group of $\h_1$.

For $H_0\neq0$ and $\h_1$ solvable (i.e., $\g$ solvable), the answer is Theorem \ref{main}, which states:
\[
\g \text{ is }\varphi\text{-symmetric}\iff
\nabla^{\h_1}R^{\h_1}=0 \ \text{ and }\ \nabla^{\h_1}\rho_0=0,
\]
where $\nabla^{\h_1}$ denotes the Levi-Civita connection on $\h_1$, with corresponding curvature $R^{\h_1}$.
In terms of the normal $j$-algebra underlying $\h_1$, it is shown that the second condition holds if and only if the mixed root spaces
$\n_{(\epsilon_i\pm\epsilon_r)/2}$ vanish for every $r>i$, where $\epsilon_i$ is the
distinguished root selected by $H_0$ (Proposition \ref{rootobstruction}). The case when the centre of $\g$ does not vanish is elementary and is recorded in Remark \ref{centre}.

The second condition does not follow from the first, and this is the point of the paper:
$\h_1$ can be locally symmetric while $\g$ fails to be $\varphi$-symmetric. The obstruction
is a single component of $\nabla\rho_0$, which we identify with a component of the curvature
derivative of the K\"ahler quotient.

For $H_0\neq0$ and $\h_1$ solvable we work with a canonical complement $\s=L(\h_1)\oplus\q$ of $\R\xi$ in $\g$.
It is an ideal of codimension one, and its group is a solvable model of the K\"ahler quotient
(Lemma \ref{Lmap}). Computing the curvature of $\s$ reduces
$\varphi$-symmetry to the two conditions above (Proposition \ref{curv}). For the converse, we
identify a complex hyperbolic factor of $\s$ with its solvable model (Proposition 
\ref{CHfactor}), and then pass from a normal $j$-algebra to a general solvable $\h_1$ by
modifying the model itself (Lemma \ref{lift}). In Theorem \ref{thm:geometric} we reinterpret Theorem \ref{main} in geometric terms: $\mathfrak g$ is $\varphi$-symmetric if and only if the normal $j$-algebra underlying $\h_1$ is locally symmetric and the distinguished root selected by $H_0$ belongs to a rank-one factor. In that case the Kähler quotient is locally isometric to $\mathbb{CH}^{n+2}\times K$, with $K$ a Hermitian symmetric space of noncompact type.

Section \ref{lowdim} contains the consequences. Combining \cite{BV1} with Proposition \ref{dim7}, every $\g$ in this class with $\dim\g\le7$ is
$\varphi$-symmetric; in dimension five this recovers \cite[Theorem
5]{LL}. The bound is sharp: in dimension nine we exhibit a structure
whose K\"ahler base is locally symmetric but which is not $\varphi$-symmetric (Example
\ref{falla-dim9}), and up to scaling its K\"ahler base is the only one in that dimension (Remark \ref{dim9unique}).
Moreover, the answer is decided not by the base alone, but by the index that $H_0$ selects: in
Example \ref{rank3}, a single normal $j$-algebra gives rise to an 11-dimensional
$\varphi$-symmetric structure and a non-$\varphi$-symmetric one, according to where $H_0$ is
placed.

\subsection*{Acknowledgements.} The authors would like to thank Giulia Dileo for suggesting the problem. The first author was partially supported by CONICET and SECYT-UNC; he would also like to thank the Abdus Salam International Centre for Theoretical Physics (ICTP), Trieste, for its hospitality during a visit in which part of this work was carried out.

\subsection*{AI use} Large language models were used during the research and preparation of this article, primarily OpenAI GPT-6 Astra and Anthropic's Claude Fable 5.1. They were used interactively to explore  research directions proposed by the authors, derive and check formulas, search for relevant literature, and test arguments for gaps. The authors checked all the mathematical arguments, computations, and references, largely wrote the article themselves, and take full responsibility for its content.

\section{Preliminaries}

\subsection{Sasakian manifolds} Let us begin by recalling some well-known definitions. Standard references for this material are \cite{Blair} and \cite{BG}.

An \textit{almost contact structure} on a $(2n+1)$-dimensional manifold $M$ is a triple $(\varphi,\eta,\xi)$ where $\xi$ is a nowhere-vanishing vector field, $\eta$ is a $1$-form and $\varphi$ is a tensor of type $(1,1)$ such that 
    \begin{equation}\label{defi variedad sasaki}
        \eta(\xi)=1, \quad \varphi^{2}=-\operatorname{I}+\xi \otimes \eta.
    \end{equation}
In this case $(M,\varphi,\eta,\xi)$ is called an \textit{almost contact manifold}. The vector field $\xi$ is known as the \textit{Reeb vector field}.

Note that condition \eqref{defi variedad sasaki} implies $\varphi(\xi)=0$ and $\eta \circ\varphi=0$. Moreover, there exists a canonical splitting of the tangent bundle of $M$ as
\[TM=\mathcal{D}\oplus \mathcal{L}\]
where $\mathcal{D}=\Ker \eta=\operatorname{Im}\varphi$ and  $\mathcal{L}$ is the line bundle generated by $\xi$.

Just as in the case of an almost complex structure, there is a notion of integrability associated with an almost contact structure: an almost contact structure $(\varphi, \eta, \xi)$ is called $\textit{normal}$ if the  Nijenhuis tensor $N_{\varphi}$ associated to the tensor $\varphi$ defined by
\begin{equation}\label{Nijenhious definicion}
    N_{\varphi}(X,Y):=\varphi^{2}[X,Y]+[\varphi X,\varphi Y]-\varphi[\varphi X,Y]-\varphi[X,\varphi Y]
\end{equation}
satisfies the condition
\begin{equation}\label{normal condicion}
    N_{\varphi}=-d\eta \otimes \xi.
\end{equation}

Condition $\eqref{normal condicion}$ is equivalent to requiring that the almost complex structure $J$ on the product $M\times \R$ defined by:
\begin{equation}\label{J en cono}
    J\left(X,f\frac{d}{d t}\right)=\left(\varphi X-f\xi,\eta(X)\frac{d}{d t}\right)
\end{equation}
be integrable, where $f$ is a smooth function on $M \times \R$ and $t$ is a coordinate on $\R$.

Now we add a Riemannian metric into the picture. A Riemannian metric $g$ on an almost contact manifold $(M,\varphi,\eta,\xi)$ is compatible with the almost contact structure if 
\[ g(\varphi X,\varphi Y)=g(X,Y) -\eta(X)\eta(Y)\]
for any vector fields $X,Y$ on $M$. In this case, $(\varphi,\eta,\xi,g)$ is called an \textit{almost contact metric structure}.

We observe that in an almost contact metric structure $(\varphi,\eta,\xi,g)$, the tensor $\varphi$ is skew-symmetric with respect to $g$ and $\xi$ is the metric dual of the $1$-form $\eta$. 
In analogy with the almost Hermitian setting, the fundamental $2$-form can be defined by
\[\Phi(X,Y):=g(X, \varphi Y).\]

\medskip

\begin{definition}
    An almost contact metric structure $(\varphi,\eta,\xi,g)$ on a $(2n+1)$-dimensional manifold $M$ is called $\textit{Sasakian}$ if it is normal and
    \begin{equation}\label{eq: d-eta}
        d\eta=2\Phi.
    \end{equation}
In this case $(M,\varphi,\eta,\xi,g)$ is called a \textit{Sasakian manifold}.
\end{definition}

Sasakian manifolds represent the most significant class of almost contact metric manifolds due to their strong connection with Kähler manifolds. In fact, the Riemannian cone of an almost contact metric manifold, equipped with an almost complex structure induced by \eqref{J en cono}, is Kähler if and only if the structure is Sasakian (see, for instance, \cite{BG}).

We point out that \eqref{eq: d-eta} implies that $\eta \wedge (d\eta)^n$ is a volume form on $M$, that is, $\eta$ is a contact form on $M$.

A Sasakian manifold $(M,\varphi,\eta,\xi,g)$ is called (locally) $\varphi$-symmetric if the geodesic reflections with respect to the integral curves of the Reeb field $\xi$ extend to local automorphisms of the Sasakian structure \cite{Ta}. Takahashi proved that this holds if and only if the local Kähler quotients of $M$ by the Reeb foliation are Hermitian locally symmetric spaces, and that when the Reeb foliation is regular, $M$ is locally $\varphi$-symmetric if and only if the quotient Kähler manifold is locally Hermitian symmetric \cite[Section 5]{Ta}. For complete simply connected manifolds, local and global $\varphi$-symmetry agree \cite{Ta}; this is the case for left-invariant structures on simply connected Lie groups, and in the cases considered in this paper the Reeb foliation is moreover regular (see \S\ref{sasakian}). The Okumura connection \cite{Ok} used in the characterisation of Sasakian $\varphi$-symmetric spaces is
\[
\widetilde\nabla_XY=\nabla_XY+\frac12\,d\eta(X,Y)\xi-\eta(X)\varphi Y+\eta(Y)\varphi X,
\]
where $\nabla$ is the Levi-Civita connection. It is a metric connection, preserves
$(\varphi,\eta,\xi)$, and has totally skew-symmetric torsion (hence, it coincides with the characteristic connection of Friedrich and Ivanov introduced in \cite{FI} for a certain family of almost contact metric structures). Takahashi's criterion \cite{Ta} gives
\begin{equation}\label{takahashi}
M \text{ is }\varphi\text{-symmetric}\iff \widetilde\nabla\widetilde R=0,
\end{equation}
where $\widetilde R$ is the curvature of $\widetilde \nabla$.

\subsection{Sasakian Lie algebras}\label{sasakian} 

We shall study left-invariant structures on Lie groups $G$. As these are uniquely
determined by their value at the group's identity, we may consider structures on Lie algebras $\mathfrak{g} = T_eG$. Let us recall that a metric Lie algebra
$(\mathfrak{g}, \langle \cdot\,,\cdot \rangle)$ is a Lie algebra $\mathfrak{g}$
equipped with a positive-definite, symmetric bilinear form
$\langle \cdot\,,\cdot \rangle$. Any such inner product corresponds uniquely to a
left-invariant Riemannian metric on each connected Lie group $G$ integrating $\mathfrak{g}$.

\begin{definition}
A \textbf{Sasakian Lie algebra} consists of a metric Lie algebra
$(\mathfrak{g}, \langle \cdot\,,\cdot \rangle)$, a vector $\xi \in \mathfrak{g}$,
its metric dual $\eta \in \mathfrak{g}^*$ and an endomorphism
$\varphi \in \mathfrak{g}^* \otimes \mathfrak{g}$ subject to:
\begin{align*}
\varphi^2 &= -\mathrm{Id} + \xi \otimes \eta, &
\langle \varphi x, \varphi y \rangle &= \langle x, y \rangle - \eta(x)\eta(y)\\
\mathrm{d}\eta(x,y) &= 2\langle x, \varphi y \rangle, &
N_\varphi &= -d\eta \otimes \xi,
\end{align*}
with $N_\varphi$ given as in \eqref{Nijenhious definicion}.
\end{definition}

As before, we set $\Phi(x,y) := \langle x, \varphi y \rangle$. Now, the normality
condition is equivalent to the pair of requirements
\[
\operatorname{ad}_\xi \circ \varphi = \varphi \circ \operatorname{ad}_\xi, \qquad
[\varphi x, \varphi y] - [x,y] = \varphi\bigl([\varphi x, y] + [x, \varphi y]\bigr)
\]
for all $x, y \in \operatorname{Ker}\eta$ (see \cite[Proposition 3.1]{AD}).

Let us name the nullspace of $\eta$ and the projection of the Lie
bracket to this nullspace by
\[
\mathfrak{h} := \operatorname{Ker}\eta, \qquad
\theta \colon \mathfrak{h} \times \mathfrak{h} \to \mathfrak{h}, \quad
\theta(\cdot\,,\cdot) := \operatorname{pr}_{\mathfrak{h}} [\cdot\,,\cdot]\big|_{\mathfrak{h}\times\mathfrak{h}},
\]
where $\operatorname{pr}_{\mathfrak{h}} \colon \mathfrak{g} \to \mathfrak{h}$ denotes
the projection coming from the decomposition
$\mathfrak{g} = \mathbb{R}\xi \oplus \mathfrak{h}$.

It was proved in \cite[Proposition~3.4]{AFV} that the centre
$\mathfrak{z}(\mathfrak{g})$ of a Sasakian Lie algebra is at most one-dimensional:
either $\mathfrak{z}(\mathfrak{g}) = 0$, or $\mathfrak{z}(\mathfrak{g}) = \mathbb{R}\xi$.
In particular, $\mathfrak{g}$ has nontrivial centre if and only if the Reeb vector
$\xi$ is central. In that case, $\theta$ is a Lie bracket on $\mathfrak{h}$, and the
Lie algebra $(\mathfrak{h} \cong \mathfrak{g}/\mathbb{R}\xi, \theta)$ inherits a
Kähler structure given by the restrictions of $\langle\cdot\,,\cdot\rangle$,
$\varphi$ and $\Phi$ to $\mathfrak{h}$. Conversely, every Kähler Lie algebra gives
rise to a Sasakian Lie algebra with nontrivial centre, namely the central extension
defined by its Kähler form \cite{AFV}.

\medskip

We recall the following notions, which will be used throughout the paper. 

\begin{definition}\label{def:locsym}
A metric Lie algebra $(\mathfrak m,\langle\cdot,\cdot\rangle)$ is \emph{locally symmetric}
if $\nabla R=0$, where $\nabla$ is the Levi-Civita connection of $\langle\cdot,\cdot\rangle$,
given on $\mathfrak m$ by the Koszul formula, and $R$ is its curvature tensor. A Kähler Lie algebra is locally symmetric if its underlying metric Lie algebra is. 
\end{definition}

The metric Lie algebra $(\mathfrak m,\langle\cdot,\cdot\rangle)$ is locally symmetric if and only if every Lie group with Lie algebra $\mathfrak m$, endowed with the corresponding left-invariant metric, is a locally symmetric Riemannian manifold; as left-invariant metrics are complete, the simply connected such group is then a
Riemannian symmetric space \cite[Ch.~XI]{KN2}. Furthermore, in the Kähler case, the simply connected Kähler Lie group is a Hermitian symmetric space.

\medskip

Now, let $\g$ be a Sasakian Lie algebra with trivial centre and $\dim\Imm\ad_\xi=2$. We use the
notation introduced in \cite[Section 4]{ACN}, where such Lie algebras are completely characterised as follows.

\begin{theorem}[{\cite[Proposition 4.2 and Theorem 4.3]{ACN}}]\label{thm:ACN}
Let $(\h_1,\theta,J,\ip{\cdot}{\cdot})$ be a K\"ahler Lie algebra with K\"ahler form
$\omega(x,y)=\ip{x}{Jy}$, and let $k\neq0$, $\gamma\in\h_1^*$ and
$H_0\in\theta(\h_1,\h_1)$ be such that
\begin{equation}\label{compat}
\mu(x)=-\tfrac12\ip{x}{JH_0},\; JH_0\perp\theta(\h_1,\h_1),\;
\theta(x,H_0)=2\mu(x)H_0,\; \omega=-\tfrac1{2k}\,d\gamma ,
\end{equation}
for all $x\in\h_1$; in particular $\h_1$ is K\"ahler-exact. Let
$\g=\R\xi\oplus\h_1\oplus\q$, with $\q=\spann\{u,v\}$, and define a bracket on $\g$ by
extending $\theta$ as follows: for $x,y\in\h_1$,
\begin{equation}\label{basic}
[x,y]=-2\omega(x,y)\xi+\theta(x,y),\qquad [\xi,x]=0,
\end{equation}
\begin{equation}\label{qbrackets}
[x,u]=\mu(x)u+\gamma(x)v,\quad [x,v]=-\gamma(x)u+\mu(x)v,
\end{equation}
\begin{equation}\label{qbrackets2}
[\xi,u]=kv, \; [\xi,v]=-ku,\; [u,v]=2\xi+H_0 .
\end{equation}
Extend the metric by declaring $\{\xi,u,v\}$ orthonormal and orthogonal to
$\h_1$, define $\varphi$ by $\varphi|_{\h_1}=J$, $\varphi\xi=0$, $\varphi u=v$, and let $\eta$ be the metric
dual of $\xi$. Then $(\g,\varphi,\eta,\xi,\ip{\cdot}{\cdot})$ is a Sasakian Lie algebra
with trivial centre and $\Imm\ad_\xi=\q$. Conversely, every Sasakian Lie algebra with
trivial centre and $\dim\Imm\ad_\xi=2$ arises in this way.
\end{theorem}

From now on $\g$ is as in Theorem \ref{thm:ACN}; thus $\q=\Imm\ad_\xi$, $v=\varphi u$ and
$J=\varphi|_{\h_1}$, and we denote by $\nabla^{\h_1}$ the Levi-Civita connection of $\h_1$.
From the Levi-Civita connection of $\g$ computed in \cite[Section 4.2]{ACN} one gets, for
$x,y,z\in\h_1$,
\begin{equation}\label{okumura}
\widetilde\nabla_xy=\nabla^{\h_1}_xy,\qquad \widetilde\nabla_\xi z=-2Jz .
\end{equation}
The last identity in \eqref{compat} reads $\gamma(\theta(x,y))=2k\,\omega(x,y)$, and it
forces $\gamma(H_0)=-2k$ whenever $H_0\neq0$. Finally, $\g$ is solvable if and only
if $H_0\neq0$ and $\h_1$ is solvable \cite[Proposition 4.5]{ACN}.

We shall repeatedly use the geometric form of \eqref{takahashi}. Let $G$ be the simply connected Lie group with Lie algebra $\mathfrak g$ and $H$ the analytic subgroup with Lie algebra $\R\xi$. In both cases treated in this paper, $H$ is closed:  if $H_0=0$, see the proof of Theorem \ref{H0zero}; if $\mathfrak g$ is solvable, $\xi$ does not belong to the nilradical of $\mathfrak g$, since $\operatorname{ad}_\xi$ is nonzero and skew-symmetric, so $H$ projects onto a closed line in the abelian group $G/N$ (with $N$ the nilradical of $G$) and is therefore closed. Consequently, the Reeb foliation is the fibration $p\colon G\to G/H$, in particular, it is regular, and $G/H$ carries a $G$-invariant Kähler structure $(\hat J,\widehat g)$ for which $p$ is a Riemannian submersion. By Takahashi’s theorem \cite{Ta} (see also \cite[Section 2]{LL} for the Lie group setting), $\mathfrak g$ is $\varphi$-symmetric if and only if $(G/H,\widehat J,\widehat g)$ is locally Hermitian symmetric. Since $G/H$ is Kähler, this is the same as local symmetry, so that
\begin{equation}\label{quotient}
\g \text{ is }\varphi\text{-symmetric}\iff (G/H,\widehat g) 
\text{ is locally  symmetric.}
\end{equation}
Left-invariant Riemannian metrics are complete, so local and global $\varphi$-symmetry agree here \cite{Ta}.

Loiudice and Lotta proved that every five-dimensional Sasakian Lie algebra with trivial
centre is $\varphi$-symmetric \cite[Theorem 5]{LL}. We study the corresponding question in higher dimensions under the assumption $\dim\Imm\ad_\xi=2$; 
the answer is given by Theorem \ref{H0zero} when $H_0=0$, and by Theorem \ref{main},
proved in Section \ref{proofmain}, when $\g$ is solvable. We close the present section with two
facts. The first one, Lemma \ref{lem:L}, will be used throughout the article; the second one is a
necessary condition for $\varphi$-symmetry.

\begin{lemma}\label{lem:L}
The linear map
\[
L\colon \mathfrak h_1\to\mathfrak g,\qquad L(x)=x-\frac{\gamma(x)}{k}\,\xi,
\]
is injective, $L(\mathfrak h_1)\cap(\mathbb R\xi\oplus\mathfrak q)=\{0\}$, and for all $x,y\in\mathfrak h_1$,
\begin{gather*}
[L(x),L(y)]=L\big(\theta(x,y)\big),\qquad [\xi,L(x)]=0,\\
[L(x),u]=\mu(x)u,\qquad [L(x),v]=\mu(x)v.
\end{gather*}
In particular, $L(\mathfrak h_1)$ is a Lie subalgebra of $\mathfrak g$ isomorphic to $(\mathfrak h_1,\theta)$.
Moreover, identifying $T_{eH}(G/H)$ with $\mathfrak h_1\oplus\mathfrak q$, one has $(dp)_eL(x)=x$ for every
$x\in\mathfrak h_1$.
\end{lemma}
 
\begin{proof}
Since $\xi$ is central in $\operatorname{Ker}\operatorname{ad}_\xi=\mathbb R\xi\oplus\mathfrak h_1$,  by \eqref{basic} and
$\gamma(\theta(x,y))=2k\,\omega(x,y)$ we have that
\[
[L(x),L(y)]=[x,y]=-2\omega(x,y)\xi+\theta(x,y)=\theta(x,y)-\frac{\gamma(\theta(x,y))}{k}\,\xi
=L\big(\theta(x,y)\big).
\]
By \eqref{qbrackets} and \eqref{qbrackets2}, $[L(x),u]=\mu(x)u+\gamma(x)v-\frac{\gamma(x)}{k}\,kv=\mu(x)u$, and similarly
$[L(x),v]=\mu(x)v$. Clearly $[\xi,L(x)]=0$.
The components of $L(x)$ are in $\mathfrak h_1\oplus\mathbb R\xi$ and its $\mathfrak h_1$-component is $x$, which gives injectivity, $L(\mathfrak h_1)\cap(\mathbb R\xi\oplus\mathfrak q)=\{0\}$ and $(dp)_eL(x)=x$,
since $(dp)_e$ vanishes on $\mathbb R\xi$.
\end{proof}

\begin{proposition}\label{necessary}
Under the assumptions above, if $\g$ is $\varphi$-symmetric, then the K\"ahler-exact Lie
algebra $\h_1$ is locally symmetric.
\end{proposition}

\begin{proof}
For $x,y,z\in\h_1$, identities \eqref{basic} and \eqref{okumura} give
\[
\begin{aligned}
\widetilde R(x,y)z
&=\widetilde\nabla_x\widetilde\nabla_yz-\widetilde\nabla_y\widetilde\nabla_xz
 -\widetilde\nabla_{[x,y]}z\\
&=\nabla^{\h_1}_x\nabla^{\h_1}_yz-\nabla^{\h_1}_y\nabla^{\h_1}_xz
 +2\omega(x,y)\widetilde\nabla_\xi z-\nabla^{\h_1}_{\theta(x,y)}z\\
&=R^{\h_1}(x,y)z-4\omega(x,y)Jz .
\end{aligned}
\]
Since $\nabla^{\h_1}\omega=0$ and $\nabla^{\h_1}J=0$, it follows that
$(\widetilde\nabla_x\widetilde R)(y,z)w=(\nabla^{\h_1}_xR^{\h_1})(y,z)w$ for
$x,y,z,w\in\h_1$. Therefore $\widetilde\nabla\widetilde R=0$ implies
$\nabla^{\h_1}R^{\h_1}=0$.
\end{proof}

\begin{remark}\label{centre}
If $\g$ has nontrivial centre then $\g=\h\oplus\R\xi$ is the central extension of a K\"ahler Lie algebra
$\h$ by the cocycle $\omega$, that is, $[x,y]=-2\omega(x,y)\xi+[x,y]_\h$. The K\"ahler
quotient is then the K\"ahler Lie group with Lie algebra $\h$. Therefore,  by \eqref{quotient}:
\[
\g\text{ is }\varphi\text{-symmetric}\iff \h\text{ is locally symmetric}.
\]
\end{remark}

\subsection{Normal \texorpdfstring{$j$}{j}-algebras and modifications}\label{sec:nja}
We recall the facts on normal $j$-algebras and their modifications that are used in Section \ref{nonzero}. All of
them are taken from \cite[Sections 5--7]{ACN}, where proofs and further references can be found.
 
\begin{definition}\label{def:nja}
A \emph{normal $j$-algebra} is a triple $(\mathfrak d,J,f)$, where $\mathfrak d$ is a completely solvable
real Lie algebra, $J$ is an integrable complex structure on $\mathfrak d$ and $f\in\mathfrak d^*$ satisfies
\[
f([Jx,Jy])=f([x,y]),\qquad f([Jx,x])>0
\]
for all $x,y\in\mathfrak d$, $x\neq0$.
\end{definition}
 
Then $\langle x,y\rangle:=f([Jx,y])$ is a $J$-invariant inner product whose Kähler form is
$\omega=-df$; thus a normal $j$-algebra is a completely solvable Kähler-exact Lie algebra. The structure of a normal $j$-algebra is completely understood, thanks to the following theorem.
 
\begin{theorem}[\cite{PS, DAtri1980}]\label{thm:nja}
Let $(\mathfrak d,J,f)$ be a normal $j$-algebra, $\mathfrak n:=[\mathfrak d,\mathfrak d]$ and
$\mathfrak a:=\mathfrak n^\perp$, so that $\mathfrak d=\mathfrak a\oplus^\perp\mathfrak n$.
\begin{enumerate}
\item[(a)] $\mathfrak a$ is abelian and $\mathfrak n=\bigoplus_{\alpha\in\Delta}\mathfrak n_\alpha$ splits into the
root spaces
\[
\mathfrak n_\alpha=\{x\in\mathfrak n:[h,x]=\alpha(h)x \text{ for all } h\in\mathfrak a\},
\]
which are mutually orthogonal and satisfy $[\mathfrak n_\alpha,\mathfrak n_\beta]\subseteq\mathfrak n_{\alpha+\beta}$.
\item[(b)] There are linearly independent roots $\epsilon_1,\dots,\epsilon_r$, $r=\dim\mathfrak a$, called
\emph{distinguished}, such that $\dim\mathfrak n_{\epsilon_k}=1$ and $J\mathfrak n_{\epsilon_k}\subset\mathfrak a$. Every
root is of the form $\epsilon_k$, $\frac12\epsilon_k$ or $\frac12(\epsilon_j\pm\epsilon_k)$ with $j<k$ (not all of
these need be roots).
\item[(c)] $J\mathfrak n_{\epsilon_k/2}=\mathfrak n_{\epsilon_k/2}$ and
$J\mathfrak n_{(\epsilon_j+\epsilon_k)/2}=\mathfrak n_{(\epsilon_j-\epsilon_k)/2}$ for $j<k$.
\item[(d)] If $j<k$ and $0\neq x\in\mathfrak n_{(\epsilon_j-\epsilon_k)/2}$, then $[x,y]\neq0$ for every non-zero
$y\in\mathfrak n_{\epsilon_k}$ and every non-zero $y\in\mathfrak n_{\epsilon_k/2}$.
\item[(e)] There are unit vectors $x_k\in\mathfrak n_{\epsilon_k}$ such that $h_k:=Jx_k$ satisfy
$\epsilon_j(h_k)=0$ for $j\neq k$ and $\epsilon_k(h_k)>0$; then $\{h_1,\dots,h_r\}$ is a basis of
$\mathfrak a$ and $\epsilon_k(h_k)\,f(x_k)=|x_k|^2=1$. Moreover, $f$ vanishes on every non-distinguished
root space.
\end{enumerate}
\end{theorem}

\begin{remark}\label{rem:solvmodel}
Let $M=G/K$ be a Hermitian symmetric space of noncompact type, with $G$ semisimple, and let $G=KAN$ be an Iwasawa decomposition \cite[Ch.~VI]{Hel}. The group $S=AN$ acts simply
transitively on $M$, so the Kähler structure of $M$ pulls back to a left-invariant Kähler
structure on $S$, and $\mathfrak s=\mathfrak a\oplus\mathfrak n$ becomes a completely solvable,
locally symmetric Kähler Lie algebra. It is well defined up to isomorphism of metric Kähler
Lie algebras, since the pair $(\mathfrak a,\mathfrak n)$ is unique up to conjugation by $K$.
Moreover, as $M$ is a bounded symmetric domain, $\mathfrak s$ is a normal $j$-algebra (\cite{PS}). We call $(\mathfrak s,J,\langle\cdot,\cdot\rangle)$ the \emph{solvable model} of $M$. For instance, it follows from \cite{He} that the solvable model of $\mathbb{C}H^{n}(-a^2)=\mathrm{SU}(n,1)/\mathrm{U}(n)$ is the algebra $\mathbb R\widehat H\oplus\mathfrak v\oplus\mathbb R\widehat Z$, $\dim\mathfrak v=2n-2$, with brackets \eqref{CHbrackets} below. 
Conversely, every locally symmetric normal $j$-algebra is the solvable model of a Hermitian symmetric space of noncompact type, see \cite{PS}.
\end{remark}

Let $(\mathfrak h,J,\langle\cdot,\cdot\rangle)$ be a Hermitian Lie algebra. We denote by
$\operatorname{Der}_u(\mathfrak h)$ the space of \emph{unitary derivations}, i.e.\ skew-symmetric derivations of
$\mathfrak h$ commuting with $J$.
 
\begin{definition}\label{def:mod}

A \emph{modification} of $\mathfrak h$ is a linear map $\phi\colon\mathfrak h\to\operatorname{Der}_u(\mathfrak h)$,
$x\mapsto\phi_x$, such that
\[
[\phi_x,\phi_y]=0,\qquad \phi_{[x,y]}=0,\qquad \phi_{\phi_xy}=0
\]
for all $x,y\in\mathfrak h$. The \emph{modified bracket} $[x,y]_\phi:=[x,y]+\phi_xy-\phi_yx$ is a Lie
bracket, and $\mathfrak h_\phi$ denotes $\mathfrak h$ with this bracket and the same $J$ and
$\langle\cdot,\cdot\rangle$.
\end{definition}
 
\begin{proposition}\label{prop:mod}
Let $\phi$ be a modification of a Hermitian Lie algebra $\mathfrak h$.
\begin{enumerate}
\item[(a)] The Levi-Civita connection and the curvature of $\mathfrak h_\phi$ are
$\nabla^\phi_x=\nabla_x+\phi_x$ and $R^\phi=R$.
\item[(b)] If $\mathfrak d$ is a normal $j$-algebra, every skew-symmetric derivation
of $\mathfrak d$ vanishes on $\mathfrak a$, and every $D\in\operatorname{Der}_u(\mathfrak d)$ vanishes on each
$\mathfrak n_{\epsilon_k}$.
\item[(c)] If $(\mathfrak d,J,f)$ is a normal $j$-algebra, every
modification $\mathfrak d_\phi$ is Kähler-exact with $\omega=-d_\phi f$. Conversely, every solvable Kähler-exact
Lie algebra is a modification of a normal $j$-algebra.
\end{enumerate}
\end{proposition}

\section{The case \texorpdfstring{$H_0=0$}{H0 = 0}}

\begin{theorem}\label{H0zero}
Under the assumptions above, suppose in addition that $H_0=0$. Then $\g$ is
$\varphi$-symmetric if and only if $\h_1$ is locally symmetric.
\end{theorem}

\begin{proof}
Since $H_0=0$,  \eqref{compat} gives $\mu=0$ and \eqref{qbrackets2} gives $[u,v]=2\xi$. Hence $\mathfrak b:=\mathbb R\xi\oplus\mathfrak q$
is a subalgebra, isomorphic to $\mathfrak{su}(2)$ if $k>0$ and to $\mathfrak{sl}(2,\mathbb R)$ if $k<0$, and by
Lemma~\ref{lem:L} the subalgebra $\tilde{\mathfrak h}_1:=L(\mathfrak h_1)\cong(\mathfrak h_1,\theta)$ satisfies
$[\mathfrak b,\tilde{\mathfrak h}_1]=0$. Since $\mathfrak b\cap\tilde{\mathfrak h}_1=\{0\}$ and
$\dim\mathfrak b+\dim\tilde{\mathfrak h}_1=\dim\mathfrak g$, we get $\mathfrak g=\mathfrak b\times\tilde{\mathfrak h}_1$.
 
Accordingly $G=B\times K$, where $B$ and $K$ are the simply connected Lie groups with Lie algebras
$\mathfrak b$ and $\tilde{\mathfrak h}_1$, and $H\subset B$. Moreover, $H$ is closed in $B$, hence in $G$: indeed, since $\operatorname{ad}_\xi|_{\mathfrak b}$ has eigenvalues $0,\pm ik$, the subgroup $H=\exp(\mathbb R\xi)$ is a maximal torus of $SU(2)$ when $k>0$, and the preimage in $\widetilde{SL}(2,\mathbb R)$ of a subgroup $SO(2)\subset SL(2,\mathbb R)$, isomorphic to $\mathbb R$, when $k<0$. Hence $G/H=(B/H)\times K$, with $G$ acting factorwise, and $B/H$ is $S^2$ or the hyperbolic plane according to the sign of $k$. Let $\widehat g_1$ and $\widehat g_2$ be the metrics induced by $\widehat g$ on the orbits
$B\cdot o\cong B/H$ and $K\cdot o\cong K$ of the origin $o=eH$. Since $G$ acts by isometries,
$\widehat g_1$ is $B$-invariant and $\widehat g_2$ is left-invariant. Under the identification
$T_o(G/H)=\mathfrak h_1\oplus\mathfrak q$, $(dp)_e$ is the identity on $\mathfrak q$ and, by Lemma~\ref{lem:L},
$(dp)_eL(x)=x$; hence $T_o(B\cdot o)=\mathfrak q$ and $T_o(K\cdot o)=\mathfrak h_1$ are orthogonal, and
$\widehat g_2(L(x),L(y))=\langle x,y\rangle$. Then $\widehat g_1\times\widehat g_2$ and $\widehat g$ are $G$-invariant and
coincide at $o$, hence they are equal.
 
Thus $(G/H,\widehat g)$ is the Riemannian product of the homogeneous surface $(B/H,\widehat g_1)$, which  has constant  curvature, and of  
$(K,\widehat g_2)$, whose metric Lie algebra is $(\mathfrak h_1 ,\langle \cdot ,\cdot \rangle)$ via $L$. A Riemannian product is locally symmetric if and only if each factor is, and \eqref{quotient} gives the claim.
\end{proof}

\section{The case \texorpdfstring{$H_0\neq0$}{H0 ≠ 0}}\label{nonzero}

Throughout this section we assume that $\g$ is solvable, equivalently, $H_0\neq0$ and $\h_1$ is solvable. In contrast with the previous section, the K\"ahler quotient $G/H$ now carries a simply transitive action of a solvable group, and the whole analysis takes place in its Lie algebra.

\subsection{A canonical complement}

\begin{lemma}\label{Lmap}
Under the assumption $H_0\neq 0$, the linear map $L\colon\h_1\to\g$ from Lemma \ref{lem:L} satisfies 
\[
L(H_0)=2\xi+H_0=[u,v] .
\]
Consequently, $\s:=L(\h_1)\oplus\q$ is an ideal of $\g$, and $\g=\R\xi\oplus\s$ as vector spaces.
\end{lemma}

\begin{proof}
The fact that $\gamma(H_0)=-2k$ (see Section \ref{sasakian}) gives immediately
\begin{equation}\label{eq: LH0} L(H_0)=H_0+2\xi=[u,v].
\end{equation}

We know from Lemma \ref{lem:L} that $L(\h_1)$ is a Lie subalgebra of $\g$. Moreover, it follows from the same lemma that $[L(\h_1),\q]\subseteq\q$, and equation \eqref{eq: LH0} gives 
$[\q,\q] \subseteq L(\h_1)$, so $\s$ is a Lie subalgebra. Since, furthermore,
$[\xi,L(\h_1)]=0$ and $[\xi,\q]\subseteq\q$, we have that $\s$ is an ideal of $\g$.

Finally, $\s\cap\R\xi=\{0\}$: if $L(x)+w=r\xi$ with $x\in\h_1$, $w\in\q$,
$r\in\R$, then
\[
x-\Big(r+\frac{\gamma(x)}{k}\Big)\xi+w=0 ,
\]
and the orthogonal decomposition $\g=\R\xi\oplus\h_1\oplus\q$ forces $x=w=0$ and then
$r=0$. Since $\dim\s=\dim\h_1+2=\dim\g-1$, we conclude $\g=\R\xi\oplus\s$.
\end{proof}

By Lemma \ref{Lmap}, $\s$ is an ideal of $\g$ and $\g=\mathbb R\xi\oplus\s$. Since $G$ is simply connected and solvable, $S$ is closed and normal, and $H\cong\mathbb R$ is closed (see \S \ref{sasakian}). As $\mathfrak g=\mathfrak s\oplus\mathbb R\xi$, we have $G=SH$, and $S\cap H$ is trivial because $H$ maps injectively onto the simply connected one-dimensional group $G/S\cong\mathbb R$. Hence $G=H\ltimes S$ and $S$ acts simply transitively on $G/H$, so that $p|_S\colon S\to G/H$ is a diffeomorphism (cf. \cite[p.~45]{LL}).
 
We transport the Kähler structure of $G/H$ to $S$ through the diffeomorphism $p|_S$. Let $P:=d(p|_S)_e=(dp)_e|_{\mathfrak s}\colon\mathfrak s\to T_o(G/H)$, which is a linear isomorphism, and endow $\mathfrak s$ with the pullback metric \[ \langle\cdot,\cdot\rangle_{\mathfrak s}:=(p|_S)^*\widehat g,\] that is, $\langle X,Y\rangle_{\mathfrak s}=\widehat g_o(PX,PY)$, and the pullback complex structure \[ J^{\mathfrak s}:=P^{-1}\circ\widehat J_o\circ P.\]
Note that $\ip{\cdot}{\cdot}_\s$ is
\emph{not} the restriction to $\s$ of the metric of $\g$: for instance
$|L(H_0)|_\g^2=4+|H_0|^2$ while $|L(H_0)|_\s^2=|H_0|^2$. This distinction is used in every connection and curvature computation below.

\begin{lemma}\label{Lisom}
For all $x,y\in\h_1$,
\[
\ip{L(x)}{L(y)}_\s=\ip{x}{y}_{\h_1},\qquad J^\s L(x)=L(Jx),\qquad J^\s u=v,
\]
and $L(\h_1)\perp\q$. In particular, $L$ is an isomorphism of metric K\"ahler Lie algebras
from $(\h_1,J,\ip{\cdot}{\cdot}_{\h_1})$ onto $(L(\h_1),J^\s,\ip{\cdot}{\cdot}_\s)$, so
$L(\h_1)$ is locally symmetric whenever $\h_1$ is.
\end{lemma}

\begin{proof}
Identify $T_{eH}(G/H)$ with $\h_1\oplus\q$. Since $L(x)-x\in\R\xi$ we have
$P(L(x))=x$ and $P|_\q=\mathrm{id}$, hence
\[
\ip{L(x)}{L(y)}_\s=\ip{P(L(x))}{P(L(y))}=\ip{x}{y}_{\h_1},
\]
and likewise $L(\h_1)\perp\q$ and
\[ J^{\mathfrak s}L(x)=P^{-1}\widehat J(x)=P^{-1}(Jx)=L(Jx),\]
while
\[ J^{\mathfrak s}u=P^{-1}\widehat J(u)=P^{-1}(\varphi u)=P^{-1}(v)=v.\]
The homomorphism property follows from Lemma \ref{lem:L}.
\end{proof}


From now on we abbreviate $\bar x:=L(x)$ for $x\in\h_1$, and we let $c\in\h_1^*$ be the
one-form
\[
c(x):=-\tfrac12\ip{x}{H_0}=\mu(Jx).
\]
The Levi-Civita connection $\nabla^\s$ of $(\s,\ip{\cdot}{\cdot}_\s)$ is computed from the
Koszul formula and Lemmas \ref{Lmap}--\ref{Lisom}. Thus, for $x,y\in\h_1$, we obtain:
\[
\begin{array}{c|ccc}
\nabla^\s & \bar y & u & v\\ \hline
\bar x & L\big(\nabla^{\h_1}_xy\big) & c(x)\,v & -c(x)\,u\\[2pt]
u & c(y)v-\mu(y)u & -\tfrac12L(JH_0) & \tfrac12L(H_0)\\[2pt]
v & -c(y)u-\mu(y)v & -\tfrac12L(H_0) & -\tfrac12L(JH_0).
\end{array}
\]
For example, $2\ip{\nabla^\s_{\bar x}u}{v}_\s=-\ip{H_0}{x}=2c(x)$.

\subsection{Curvature and the two-form \texorpdfstring{$\rho_0$}{ρ0}}

We now introduce the two-form on $\h_1$ that records the component of the bracket along $H_0$:
\[
\rho_0(x,y):=\ip{\theta(x,y)}{H_0},\qquad x,y\in\h_1 .
\]

\begin{proposition}\label{curv}
For all $x,y,z,w\in\h_1$,
\[
R^\s(\bar x,\bar y)\bar z=L\big(R^{\h_1}(x,y)z\big),\qquad
R^\s(\bar x,\bar y)u=\tfrac12\rho_0(x,y)\,v,\qquad
R^\s(\bar x,\bar y)v=-\tfrac12\rho_0(x,y)\,u ,
\]
and consequently
\[
(\nabla^\s_{\bar z}R^\s)(\bar x,\bar y)\bar w=L\big((\nabla^{\h_1}_zR^{\h_1})(x,y)w\big),
\qquad
(\nabla^\s_{\bar z}R^\s)(\bar x,\bar y)u=\tfrac12\big(\nabla^{\h_1}_z\rho_0\big)(x,y)\,v .
\]
In particular, if $\nabla^\s R^\s=0$ then $\h_1$ is locally symmetric and
$\nabla^{\h_1}\rho_0=0$.
\end{proposition}

\begin{proof}
The first identity is immediate from the connection table and
$[\bar x,\bar y]=L(\theta(x,y))$. For the second,
\[
\begin{aligned}
R^\s(\bar x,\bar y)u&=\nabla^\s_{\bar x}\big(c(y)v\big)-\nabla^\s_{\bar y}\big(c(x)v\big)
-\nabla^\s_{L(\theta(x,y))}u\\
&=-c(y)c(x)u+c(x)c(y)u-c\big(\theta(x,y)\big)v=\tfrac12\rho_0(x,y)\,v ,
\end{aligned}
\]
and analogously $R^\s(\bar x,\bar y)v=-\tfrac12\rho_0(x,y)u$. Differentiating the first
identity gives the third one. Finally, differentiating the second,
\[
\begin{aligned}
(\nabla^\s_{\bar z}R^\s)(\bar x,\bar y)u
&=\tfrac12\rho_0(x,y)\nabla^\s_{\bar z}v-\tfrac12\rho_0\big(\nabla^{\h_1}_zx,y\big)v
 -\tfrac12\rho_0\big(x,\nabla^{\h_1}_zy\big)v-c(z)R^\s(\bar x,\bar y)v\\
&=-\tfrac12\Big\{\rho_0\big(\nabla^{\h_1}_zx,y\big)+\rho_0\big(x,\nabla^{\h_1}_zy\big)\Big\}v
=\tfrac12\big(\nabla^{\h_1}_z\rho_0\big)(x,y)\,v ,
\end{aligned}
\]
the two terms in $c(z)$ cancelling.
\end{proof}

\subsection{Metric lemmas}

\begin{lemma}\label{deriv}
Let $(\m,\ip{\cdot}{\cdot})$ be a metric Lie algebra, with Levi-Civita connection $\nabla$
and curvature tensor $R$, and let $D$ be a skew-symmetric derivation of $\m$. Then, for all
$x,y,z\in\m$,
\begin{equation}\label{Dnabla}
D\big(\nabla_xy\big)=\nabla_{Dx}y+\nabla_x(Dy),
\end{equation}
\begin{equation}\label{DR}
D\big(R(x,y)z\big)=R(Dx,y)z+R(x,Dy)z+R(x,y)Dz .
\end{equation}
\end{lemma}

\begin{proof}
Write $A:=\nabla_{Dx}y+\nabla_x(Dy)$. By the Koszul formula,
\[
-2\ip{\nabla_xy}{Dz}=-\ip{[x,y]}{Dz}+\ip{[y,Dz]}{x}-\ip{[Dz,x]}{y},
\]
and, again by Koszul,
\[
\begin{aligned}
2\ip{A}{z}&=\ip{[Dx,y]}{z}-\ip{[y,z]}{Dx}+\ip{[z,Dx]}{y}
+\ip{[x,Dy]}{z}-\ip{[Dy,z]}{x}+\ip{[z,x]}{Dy}\\
&=\ip{D[x,y]}{z}+\ip{D[y,z]-[Dy,z]}{x}+\ip{[z,Dx]-D[z,x]}{y}\\
&=-\ip{[x,y]}{Dz}+\ip{[y,Dz]}{x}-\ip{[Dz,x]}{y},
\end{aligned}
\]
where the derivation identity was used in the second equality and skew-symmetry of $D$ in
the third. Therefore $-2\ip{\nabla_xy}{Dz}=2\ip{A}{z}$ for all $z\in\m$, that is,
$D(\nabla_xy)=A$, which is \eqref{Dnabla}.

For \eqref{DR}, apply $D$ to $R(x,y)z=\nabla_x\nabla_yz-\nabla_y\nabla_xz-\nabla_{[x,y]}z$
and use \eqref{Dnabla} twice in each of the first two terms:
\[
D\big(\nabla_x\nabla_yz\big)=\nabla_{Dx}\nabla_yz+\nabla_x\nabla_{Dy}z+\nabla_x\nabla_y(Dz),
\]
and similarly with $x$ and $y$ interchanged, while
\[
D\big(\nabla_{[x,y]}z\big)=\nabla_{D[x,y]}z+\nabla_{[x,y]}(Dz)
=\nabla_{[Dx,y]}z+\nabla_{[x,Dy]}z+\nabla_{[x,y]}(Dz),
\]
because $D$ is a derivation. Collecting the nine resulting terms into those in $Dx$, those
in $Dy$ and those in $Dz$ gives
\[
\big(\nabla_{Dx}\nabla_yz-\nabla_y\nabla_{Dx}z-\nabla_{[Dx,y]}z\big)
+\big(\nabla_x\nabla_{Dy}z-\nabla_{Dy}\nabla_xz-\nabla_{[x,Dy]}z\big)+R(x,y)Dz,
\]
that is, $R(Dx,y)z+R(x,Dy)z+R(x,y)Dz$.
\end{proof}

\begin{corollary}\label{Dacts}
Let $(\m,\ip{\cdot}{\cdot})$ be a metric Lie algebra and let $D$ be a skew-symmetric
derivation of $\m$. Then $D\cdot R=0$, and $D\cdot\sigma=0$ for every bilinear form
$\sigma(x,y)=\ip{[x,y]}{H}$ with $DH=0$. Here $\cdot$ denotes the natural action of $D$ as
a derivation on tensors,
\begin{align*}
    (D\cdot R)(x,y)z&=D\big(R(x,y)z\big)-R(Dx,y)z-R(x,Dy)z-R(x,y)Dz,\\
    (D\cdot\sigma)(x,y)&=-\sigma(Dx,y)-\sigma(x,Dy).
\end{align*}
\end{corollary}

\begin{proof}
The first assertion is \eqref{DR}. For the second, using that $D$ is a derivation and
skew-symmetric,

\[
(D\cdot\sigma)(x,y)=-\ip{[Dx,y]+[x,Dy]}{H}=-\ip{D[x,y]}{H}=\ip{[x,y]}{DH}=0 .
\qedhere
\]
\end{proof}

We shall also use the following elementary fact. If a metric Lie algebra
$(\m,\ip{\cdot}{\cdot})$ decomposes as an orthogonal direct sum $\m=\m_1\oplus^\perp\m_2$
of commuting ideals, the Koszul formula gives $\nabla^\m_{\m_j}\m_j\subseteq\m_j$ and
$\nabla^\m_{\m_1}\m_2=\nabla^\m_{\m_2}\m_1=0$, so that
\[
\nabla^\m|_{\m_j\times\m_j}=\nabla^{\m_j} ,
\qquad
R^\m=R^{\m_1}\oplus R^{\m_2} ,
\qquad
\nabla^\m R^\m=\nabla^{\m_1}R^{\m_1}\oplus\nabla^{\m_2}R^{\m_2} ,
\]
all mixed components being zero; in particular $\m$ is locally symmetric if and only if
$\m_1$ and $\m_2$ are.

\subsection{Normal \texorpdfstring{$j$}{j}-algebras and the mixed root obstruction}\label{mixed}

We now exploit the solvability assumption and use Proposition~\ref{prop:mod}, which allows us to consider $\h_1$ as a modification of a normal $j$-algebra $\mathfrak d=\mathfrak a\oplus^\perp\mathfrak n$, where we are using the notation of Theorem~\ref{thm:nja}. By \cite[Proposition~7.5]{ACN} there exists an index $i$ such
that
\[
H_0=-\epsilon_i(h_i)\,x_i,
\]
and we write $a:=\epsilon_i(h_i)=|H_0|>0$.

Moreover, $J$ preserves $\mathfrak n_{\epsilon_j/2}$
for every $j$, and interchanges
$\mathfrak n_{(\epsilon_j-\epsilon_k)/2}$ and
$\mathfrak n_{(\epsilon_j+\epsilon_k)/2}$
for every pair of indices $j<k$.

We first state a useful lemma that shows that certain root spaces vanish. This result shows that the $\eta$-Einstein hypothesis is not necessary in \cite[Proposition 7.7]{ACN}.

\begin{lemma}\label{smaller}
If $r<i$, then $\n_{(\epsilon_r-\epsilon_i)/2}=\n_{(\epsilon_r+\epsilon_i)/2}=\{0\}$.
\end{lemma}

\begin{proof} Suppose that $0\neq w\in\n_{(\epsilon_r-\epsilon_i)/2}$. Since $JH_0=-a\,h_i\in \aaa$ and $w\in\n=\aaa^\perp$, we have $\mu(w)=-\tfrac12\ip{w}{JH_0}=0$, so the
compatibility condition \eqref{compat} gives
\[
\theta(w,H_0)=2\mu(w)H_0=0 .
\]
On the other hand, Theorem~\ref{thm:nja}, applied to the pair $(r,i)$ with $w\in\n_{(\epsilon_r-\epsilon_i)/2}$ and $x_i\in\n_{\epsilon_i}$, implies
$\theta(w,x_i)=[w,x_i]_{\mathfrak d}\neq0$ (since a modification vanishes on the commutator ideal) and hence $\theta(w,H_0)\neq0$, a contradiction. Hence, $\n_{(\epsilon_r-\epsilon_i)/2}=0$, and $\n_{(\epsilon_r+\epsilon_i)/2}=J\n_{(\epsilon_r-\epsilon_i)/2}=0$.
\end{proof}

Recall also from Theorem~\ref{thm:nja} that $f(x_i)=1/a$. By \cite[p.~65]{DAtri1980}, the Levi-Civita connection of the normal $j$-algebra $\mathfrak d$ satisfies
\begin{equation}\label{nablaa}
    \nabla^{\mathfrak d}_{h}=0
    \qquad\text{for every }h\in\mathfrak a.
\end{equation}

\begin{lemma}\label{modbase}
Let $\phi$ be a modification of a normal $j$-algebra $\dd$ and $\h_1=\dd_\phi$. Then $\h_1$
is locally symmetric if and only if $\dd$ is, the form $\rho_0$ is the same for the two
brackets, and $\nabla^{\h_1}\rho_0=\nabla^{\dd}\rho_0$.
\end{lemma}

\begin{proof}

 By Definition~\ref{def:mod} each $\phi_x$ is a skew-symmetric derivation of $\mathfrak d$ commuting with $J$,
and by Proposition~\ref{prop:mod} the Levi-Civita connection $\nabla^\phi_x$ of $\mathfrak h_1=\mathfrak d_\phi$
and its curvature $R^\phi$ are given by
\[
\nabla^\phi_x=\nabla_x+\phi_x,\qquad R^\phi=R^{\dd}.
\]
Applying this to $R^\dd$ and using Corollary \ref{Dacts},
\[
\nabla^\phi_xR^\phi=\nabla_xR^{\dd}+\phi_x\cdot R^{\dd}=\nabla_xR^{\dd},
\]
whence $\nabla^{\h_1}R^{\h_1}=0$ if and only if $\nabla^{\dd}R^{\dd}=0$.

For the second claim, by Proposition~\ref{prop:mod} and $H_0\in\mathfrak n_{\epsilon_i}$ we have $\phi_xH_0=0$, so
\[
\ip{\phi_xy}{H_0}=-\ip{y}{\phi_xH_0}=0 ,
\]
and $[x,y]_\phi=[x,y]+\phi_xy-\phi_yx$ gives $\ip{[x,y]_\phi}{H_0}=\ip{[x,y]}{H_0}$. Thus
$\rho_0$ is the same for $\dd$ and for $\h_1$. Corollary \ref{Dacts} applied to $\sigma=\rho_0$ with $H=H_0$ yields
\[
\begin{aligned}
(\nabla^\phi_x\rho_0)(y,z)&=-\rho_0(\nabla^\phi_xy,z)-\rho_0(y,\nabla^\phi_xz)\\
&=\big\{-\rho_0(\nabla^{\dd}_xy,z)-\rho_0(y,\nabla^{\dd}_xz)\big\}
+\big\{-\rho_0(\phi_xy,z)-\rho_0(y,\phi_xz)\big\}\\
&=(\nabla^{\dd}_x\rho_0)(y,z)+(\phi_x\cdot\rho_0)(y,z)=(\nabla^{\dd}_x\rho_0)(y,z). 
\end{aligned}
\]
\end{proof}

By Lemma \ref{modbase} both the local symmetry of $\h_1$ and the parallelism of $\rho_0$
may be tested on the underlying normal $j$-algebra. \emph{Throughout the rest of this subsection, we therefore
assume that $\h_1$ itself is a normal $j$-algebra}, with bracket $\theta$, and we keep the
notation
\[
H_0=-a\,x_i ,\qquad h_i=Jx_i ,\qquad |x_i|=1 ,\qquad a=\epsilon_i(h_i)>0 .
\]

The following observation will be very useful in computations below: the relation $\omega=-df$ gives $f(\theta(x,y))=\omega(x,y)$, so if $\theta(x,y)=c\,x_i$ then
\[
c\,f(x_i)=f\big(\theta(x,y)\big)=\omega(x,y) ,
\]
that is, $c=a\,\omega(x,y)$ because $f(x_i)=1/a$; pairing with $H_0=-a\,x_i$ then computes
$\rho_0$ as well. Thus, for all $x,y\in\h_1$,
\begin{equation}\label{ftrick}
\theta(x,y)\in\n_{\epsilon_i}
\quad\Longrightarrow\quad
\left\{
\begin{aligned}
\theta(x,y)&=a\,\omega(x,y)\,x_i ,\\[3pt]
\rho_0(x,y)&=-a^2\,\omega(x,y) .
\end{aligned}
\right.
\end{equation}

We now split off the part of $\h_1$ singled out by the distinguished root $\epsilon_i$: we
set
\[
\dd_i:=\spann\{h_i,x_i\}\oplus^\perp\n_{\epsilon_i/2} ,
\qquad
\kk_i:=\dd_i^\perp .
\]

\begin{lemma}\label{lem:rho0-components}
The decomposition $\h_1=\dd_i\oplus\kk_i$ is orthogonal, both summands are
$J$-invariant, $H_0, JH_0 \in\dd_i$ and $\dd_i$ is a Lie subalgebra. Moreover
\begin{equation}\label{rhodi}
\rho_0|_{\dd_i}=-a^2\,\omega|_{\dd_i} ,\qquad \rho_0(\dd_i,\kk_i)=0 .
\end{equation}
\end{lemma}

\begin{proof}
Since $Jx_i=h_i$, $Jh_i=-x_i$ and $J$ preserves $\mathfrak n_{\epsilon_i/2}$,
the subspace $\mathfrak d_i$ is $J$-invariant, and so is $\mathfrak k_i$
because $J$ is orthogonal. Also, $H_0=-a\,x_i$ and $JH_0=-a\,h_i$ lie in $\mathfrak d_i$. Now, the brackets
\[
\theta(h_i,x_i)=a\,x_i ,
\qquad
\theta(h_i,x)=\tfrac a2\,x ,
\qquad
\theta(x_i,x)=0
\qquad (x\in\n_{\epsilon_i/2}) ,
\]
together with $\theta(\n_{\epsilon_i/2},\n_{\epsilon_i/2})\subseteq\n_{\epsilon_i}$, give
$\theta(\dd_i,\dd_i)\subseteq\dd_i$, that is, $\dd_i$ is a Lie subalgebra.
 
Call a vector of $\mathfrak h_1$ \emph{homogeneous} if it belongs to
$\mathfrak a$ or to a root space $\mathfrak n_\alpha$. Both $\mathfrak d_i$ and
$\mathfrak k_i$ are spanned by homogeneous vectors:
\[
\mathfrak d_i=\mathbb R h_i\oplus\mathfrak n_{\epsilon_i}\oplus\mathfrak n_{\epsilon_i/2},
\qquad
\mathfrak k_i=\bigoplus_{j\neq i}\mathbb R h_j\;\oplus\bigoplus_{\alpha\neq\epsilon_i,\,\epsilon_i/2}\mathfrak n_\alpha .
\]
Since both sides of \eqref{rhodi} are bilinear, it suffices to check
\eqref{rhodi} on pairs of homogeneous vectors.
 
\smallskip
Note that for homogeneous $x,y$ the bracket
$\theta(x,y)$ lies in a single root space or vanishes, because
$\theta(\mathfrak a,\mathfrak a)=0$,
$\theta(\mathfrak a,\mathfrak n_\beta)\subseteq\mathfrak n_\beta$ and
$\theta(\mathfrak n_\alpha,\mathfrak n_\beta)\subseteq\mathfrak n_{\alpha+\beta}$,
the latter being zero when $\alpha+\beta$ is not a root. As the root spaces
are mutually orthogonal and $\rho_0(x,y)=-a\langle\theta(x,y),x_i\rangle$,
we get, for all homogeneous $x,y$,
\begin{equation}\label{eq:hom}
\rho_0(x,y)=
\begin{cases}
-a^2\,\omega(x,y) & \text{if }\theta(x,y)\in\mathfrak n_{\epsilon_i},\\[2pt]
0 & \text{otherwise},
\end{cases}
\end{equation}
the first case being \eqref{ftrick}.
 
\smallskip

We prove next that $\rho_0(\mathfrak d_i,\mathfrak k_i)=0$. Let
$x\in\mathfrak d_i$ and $y\in\mathfrak k_i$ be homogeneous. By
\eqref{eq:hom}, either $\rho_0(x,y)=0$ or
$\rho_0(x,y)=-a^2\omega(x,y)=-a^2\langle x,Jy\rangle$, and the latter
vanishes because $Jy\in\mathfrak k_i\perp\mathfrak d_i$.
 
\smallskip

Finally, we prove that $\rho_0|_{\mathfrak d_i}=-a^2\omega|_{\mathfrak d_i}$.
Let $x,y$ be homogeneous vectors of $\mathfrak d_i$. If
$\theta(x,y)\in\mathfrak n_{\epsilon_i}$, the claim is \eqref{eq:hom}. This
covers the pairs $(h_i,x_i)$, for which $\theta(h_i,x_i)=a\,x_i$, and the
pairs of vectors of $\mathfrak n_{\epsilon_i/2}$, for which
$\theta(\mathfrak n_{\epsilon_i/2},\mathfrak n_{\epsilon_i/2})\subseteq\mathfrak n_{\epsilon_i}$.
In the remaining cases, $\rho_0(x,y)=0$ by \eqref{eq:hom}, and we must check
that $\omega(x,y)=0$ as well:
\begin{itemize}
\item pairs $(x,x)$ with $x\in\{h_i,x_i\}$: trivially $\omega(x,x)=0$;
\item $(h_i,y)$ with $y\in\mathfrak n_{\epsilon_i/2}$: here
$\theta(h_i,y)=\tfrac a2\,y\in\mathfrak n_{\epsilon_i/2}$, and
$\omega(h_i,y)=\langle h_i,Jy\rangle=0$ because
$Jy\in\mathfrak n_{\epsilon_i/2}\perp\mathfrak a$;
\item $(x_i,y)$ with $y\in\mathfrak n_{\epsilon_i/2}$: here
$\theta(x_i,y)=0$, since $\epsilon_i+\tfrac12\epsilon_i$ is not a root, and
$\omega(x_i,y)=\langle x_i,Jy\rangle=0$ because
$Jy\in\mathfrak n_{\epsilon_i/2}\perp\mathfrak n_{\epsilon_i}$.
\end{itemize}
Hence, $\rho_0$ and $-a^2\omega$ agree on all pairs of homogeneous vectors of
$\mathfrak d_i$, and therefore on $\mathfrak d_i$. 
\end{proof}

\begin{lemma}\label{lem:rho0-commuting-ideals}
Assume that $\n_{(\epsilon_i-\epsilon_r)/2}=\n_{(\epsilon_i+\epsilon_r)/2}=\{0\}$ for every
$r>i$. Then $\dd_i$ and $\kk_i$ are commuting ideals of $\h_1$,
\begin{equation}\label{rhosplit}
\rho_0=-a^2\,\omega|_{\dd_i}\oplus 0 ,
\end{equation}
and consequently $\nabla^{\h_1}\rho_0=0$.
\end{lemma}

\begin{proof}
By Lemma \ref{smaller} and the hypothesis  all the mixed root spaces involving $\epsilon_i$ vanish, so every
root occurring in $\kk_i$ vanishes on $h_i$:
\[
\kk_i=\spann\{h_j:\ j\neq i\}
\oplus^\perp
\bigoplus_{\substack{\alpha\in\Delta\\ \alpha(h_i)=0}}\n_\alpha .
\]
By Lemma \ref{lem:rho0-components}, $\dd_i$ is a Lie subalgebra. Now, the bracket of two root vectors has
root the sum of their roots, whenever that sum is a root; if both roots vanish on $h_i$, so
does their sum. Together with the action of $\aaa$ on each root space, this gives
$\theta(\kk_i,\kk_i)\subseteq\kk_i$, so that $\kk_i$ is also a Lie subalgebra.

Let us now show that $\dd_i$ and $\kk_i$ commute. First,
\[
\theta(h_j,\dd_i)=0 \quad (j\neq i) ,
\qquad
\theta(h_i,\kk_i)=0 ,
\]
because every root occurring in $\dd_i$ vanishes on $h_j$ for $j\neq i$, and every root
occurring in $\kk_i$ vanishes on $h_i$. Second, adding $\epsilon_i$ or $\frac12 \epsilon_i$ to a
root occurring in $\kk_i$ never produces a root, the only exception being
$\tfrac12\epsilon_i+\tfrac12\epsilon_j=\tfrac12(\epsilon_i+\epsilon_j)$ with $j\neq i$,
whose root space vanishes by  Lemma \ref{smaller} and  hypothesis; hence
\[
\theta(x_i,\kk_i)=0 ,
\qquad
\theta(\n_{\epsilon_i/2},\kk_i)=0 .
\]
Therefore $\theta(\dd_i,\kk_i)=0$, and as a consequence, both $\dd_i$ and $\kk_i$ are ideals.

Since $\theta(\kk_i,\kk_i) \subseteq\kk_i\perp H_0$, we get
$\rho_0|_{\kk_i}=0$, which together with \eqref{rhodi} is exactly
\eqref{rhosplit}. Finally, $\dd_i$ and $\kk_i$ being commuting orthogonal ideals, the
Levi-Civita connection of $\h_1$
preserves both factors; as $\omega|_{\dd_i}$ is the K\"ahler form of $\dd_i$ it is parallel,
and so is the zero form on $\kk_i$. Hence $\nabla^{\h_1}\rho_0=0$.
\end{proof}

\begin{proposition}\label{rootobstruction}
Let $\h_1$ be a normal $j$-algebra and suppose that $H_0=-a\,x_i$ with $a=\epsilon_i(h_i)>0$,
$|x_i|=1$, $h_i=Jx_i$. Then
\[
\nabla^{\h_1}\rho_0=0
\iff
\n_{(\epsilon_i-\epsilon_r)/2}=\n_{(\epsilon_i+\epsilon_r)/2}=\{0\}
\quad\text{for every } r>i .
\]
\end{proposition}

\begin{proof}
($\Leftarrow$) This is Lemma \ref{lem:rho0-commuting-ideals}.

\medskip
($\Rightarrow$) Assume that $\nabla^{\h_1}\rho_0=0$, fix $r>i$ and suppose that
$0\neq w\in\n_\alpha$ with $\alpha=\tfrac12(\epsilon_i-\epsilon_r)$. We compute
$(\nabla^{\h_1}_w\rho_0)(h_r,Jw)$.

Using that $\nabla^{\h_1}$ is torsion-free and \eqref{nablaa} we obtain
\[
\nabla^{\h_1}_w h_r=\nabla^{\h_1}_{h_r}w+\theta(w,h_r)=-\alpha(h_r)\,w .
\]
Next, let us see that $\rho_0(h_r,y)=0$ for every $y\in\h_1$. Indeed,  this is clear for $y\in\aaa$, because $\aaa$ is
abelian, and if $y\in\n_\beta$ then
\[
\rho_0(h_r,y)=\ip{\theta(h_r,y)}{H_0}=\beta(h_r)\,\ip{y}{H_0}=0 ,
\]
since $\ip{y}{H_0}\neq0$ forces $\beta=\epsilon_i$, while $\epsilon_i(h_r)=0$ for $r\neq i$.

Thus, we have 
\begin{align*}\label{rhomixed}
(\nabla^{\h_1}_w\rho_0)(h_r,Jw)
&=-\rho_0\big(\nabla^{\h_1}_wh_r,Jw\big)-\rho_0\big(h_r,\nabla^{\h_1}_wJw\big)\nonumber\\[4pt]
&=\alpha(h_r)\,\rho_0(w,Jw)\\[4pt]
&=\alpha(h_r)\,a^2|w|^2,\nonumber
\end{align*}
where in the last step we used \eqref{ftrick}, since $\theta(w,Jw)\in\n_{\epsilon_i}$ and
$\omega(w,Jw)=-|w|^2$. The left-hand side vanishes by assumption, while
$\alpha(h_r)=-\tfrac12\epsilon_r(h_r)\neq0$ and $a>0$; this forces $w=0$, a
contradiction. Finally, $J$ interchanges $\n_{(\epsilon_i-\epsilon_r)/2}$ and
$\n_{(\epsilon_i+\epsilon_r)/2}$, so both spaces vanish.
\end{proof}

\subsection{The complex hyperbolic factor}

\begin{proposition}\label{CHfactor}
Assume that $\h_1$ is a locally symmetric normal $j$-algebra, that $H_0=-a\,x_i$ with
$a=\epsilon_i(h_i)=|H_0|$, and that
$\n_{(\epsilon_i-\epsilon_r)/2}=\n_{(\epsilon_i+\epsilon_r)/2}=\{0\}$ for every $r > i$. Then
the metric Lie algebra $\s$ admits an orthogonal decomposition into commuting
$J^\s$-invariant ideals
\[
\s=\s_i\oplus^\perp\widetilde\kk_i,\qquad \s_i:=L(\dd_i)\oplus\q,\qquad
\widetilde\kk_i:=L(\kk_i),
\]
where $\widetilde\kk_i$ is isometric to $\kk_i$ and $\s_i$ is the solvable model of
$\CH^{n+2}(-a^2)$, with $2n=\dim\n_{\epsilon_i/2}$. In particular $\nabla^\s R^\s=0$.
\end{proposition}

\begin{proof}
By Lemma  \ref{lem:rho0-commuting-ideals}, $\h_1=\dd_i\oplus^\perp\kk_i$. Since $JH_0\in\dd_i$
we have $JH_0\perp\kk_i$, therefore $\mu(x)=-\tfrac12\ip{x}{JH_0}=0$ for every $x\in\kk_i$,
and Lemma \ref{lem:L} gives
\[
[\widetilde\kk_i,\q]=[L(\kk_i),\q]=\mu(\kk_i)\q=0 .
\]
Moreover
\[
[\s_i,\widetilde\kk_i]=[L(\dd_i),L(\kk_i)]+[\q,\widetilde\kk_i]
=L\big(\theta(\dd_i,\kk_i)\big)=0 ,
\]
due to Lemma \ref{lem:rho0-commuting-ideals}.
Now, it follows from Lemma \ref{Lmap} that  $[\s_i,\s_i]\subseteq\s_i$  and
$[\widetilde\kk_i,\widetilde\kk_i]\subseteq\widetilde\kk_i$. Since $\h_1=\dd_i\oplus\kk_i$
and $\s=L(\h_1)\oplus\q$, we obtain $\s=\s_i\oplus\widetilde\kk_i$ with $\s_i$ and
$\widetilde\kk_i$ ideals of $\s$, orthogonal by Lemma \ref{Lisom} and $J^\s$-invariant
because $\dd_i,\kk_i$ are $J$-invariant and $J^\s\q=\q$.

Put
\[
V:=\n_{\epsilon_i/2},\qquad A:=L(H_0),\qquad B:=L(JH_0) .
\]
Since $|A|_\s=|B|_\s=a$ by Lemma \ref{Lisom}, we define the unit vectors
\[
\widehat H:=-\frac1aB,\qquad \widehat Z:=\frac1aA ,
\]
and the subspace $\mathfrak v:=L(V)\oplus\q$, so that
$J^\s\widehat H=\widehat Z$ and $J^\s\mathfrak v=\mathfrak v$. The bracket relations give
\[
\begin{aligned}
\relax[\widehat H,\widehat Z]&=-\frac1{a^2}[B,A]=-\frac1{a^2}L\big(\theta(JH_0,H_0)\big)\\
&=-\frac1{a^2}L\big(2\mu(JH_0)H_0\big)=-\frac{2\mu(JH_0)}{a^2}\,L(H_0)=a\widehat Z .
\end{aligned}
\]
For $w\in\mathfrak v$:
\begin{itemize}
\item if $w=L(x)$ with $x\in V$, then
\[
[\widehat H,L(x)]=-\frac1aL\big(\theta(JH_0,x)\big)=-\frac1aL\Big(-\frac{a^2}{2}x\Big)
=\frac a2L(x),
\]

\item if $w=u$, then $[\widehat H,u]=-\tfrac1a\mu(JH_0)u=\tfrac a2u$, and similarly for $v$.
\end{itemize}
Therefore $[\widehat H,w]=\tfrac a2w$ for every $w\in\mathfrak v$, and clearly
$[\widehat Z,\mathfrak v]=0$. On the other hand, for $w,w'\in\mathfrak v$:
\begin{itemize}
\item if $w=L(x)$, $w'=L(y)$ with $x,y\in V$, then, by \eqref{ftrick},
\[
[w,w']=L\big(\theta(x,y)\big)=L\big(a\omega(x,y)x_i\big)=-a\omega(x,y)\widehat Z
=a\ip{J^\s w}{w'}_\s\widehat Z ;
\]
\item if $w\in L(V)$ and $w'\in\q$, then $[w,w']=\mu(V)\q=0$, and also
$\ip{J^\s w}{w'}_\s=0$ since $L(\h_1)\perp\q$;
\item if $w,w'\in\q$, then $[u,v]=L(H_0)=a\widehat Z=a\ip{J^\s u}{v}_\s\widehat Z$.
\end{itemize}
In summary $\s_i=\R\widehat H\oplus^\perp\mathfrak v\oplus^\perp\R\widehat Z$ with 
\begin{equation}\label{CHbrackets}
[\widehat H,\widehat Z]=a\widehat Z,\qquad
[\widehat H,w]=\frac a2\,w,\qquad
[\widehat Z,w]=0,\qquad
[w,w']=a\ip{J^\s w}{w'}_\s\,\widehat Z,
\end{equation}
for $w,w'\in\mathfrak v$; in particular, $\dim\s_i=2n+4$.


Therefore, $(\s_i,\ip{\cdot}{\cdot}_\s)$ is the solvable model of the complex hyperbolic space of
real dimension $2n+4$; equivalently, \eqref{CHbrackets} are the standard brackets of $\CH^{n+2}$ in its solvable presentation (see, for instance, \cite{He}). Thus, $\s_i$
is locally symmetric, i.e., $\nabla^{\s_i}R^{\s_i}=0$. Its simply connected Lie group is therefore a Riemannian symmetric space.

Furthermore, by Lemmas \ref{Lmap} and \ref{Lisom}, $L$ is an injective isometric homomorphism, so $\widetilde\kk_i$ is
isometric to $\kk_i$ and
\[
\nabla^{\widetilde\kk_i}_{L(x)}L(y)=L\big(\nabla^{\kk_i}_xy\big),\qquad
R^{\widetilde\kk_i}\big(L(x),L(y)\big)L(z)=L\big(R^{\kk_i}(x,y)z\big).
\]
Since $\h_1$ is locally symmetric and $\h_1=\dd_i\oplus^\perp\kk_i$ is an orthogonal sum of
commuting ideals, $\kk_i$ is locally symmetric, hence so is $\widetilde\kk_i$. Applying the
same splitting of $\nabla R$ to $\s=\s_i\oplus^\perp\widetilde\kk_i$,
\[
\nabla^\s R^\s=\big(\nabla^{\s_i}R^{\s_i}\big)\oplus
\big(\nabla^{\widetilde\kk_i}R^{\widetilde\kk_i}\big)=0.
\]


It remains to identify the constant of the complex hyperbolic space, and for that it suffices to evaluate the curvature on one holomorphic plane. We take the plane spanned by the orthonormal pair $\widehat H$ and
$J^\s\widehat H=\widehat Z$. Using the Koszul formula we obtain 
\[
\nabla^\s_{\widehat H}\widehat H=0, \quad \nabla^\s_{\widehat H}\widehat Z=0, \quad \nabla^\s_{\widehat Z}\widehat Z=a\widehat H.
\]
Consequently
\[
R^{\s_i}(\widehat H,\widehat Z)\widehat Z
=\nabla^\s_{\widehat H}\big(a\widehat H\big)-\nabla^\s_{\widehat Z}\big(0\big)
-\nabla^\s_{a\widehat Z}\widehat Z
=-a^2\,\widehat H ,
\]
and therefore
\[
K^{\s_i}\big(\widehat H,J^\s\widehat H\big)
=\ip{R^{\s_i}(\widehat H,\widehat Z)\widehat Z}{\widehat H}_\s=-a^2 .
\]
Thus $\s_i$ is the solvable model of $\CH^{n+2}(-a^2)$.
\end{proof}

\subsection{Lifting modifications to the quotient model}

Proposition \ref{CHfactor} assumes that $\h_1$ is a normal $j$-algebra. To reach the general
solvable case we must transfer the modification from $\h_1$ to the quotient model $\s$.

Let $\h_1=\dd_\phi$ with $(\dd,J,f)$ a normal $j$-algebra, and recall from Section \ref{mixed} that
$H_0=-a\,x_i\in\n_{\epsilon_i}$. We first show that the data $k$, $\gamma$ and $H_0$ of $\g$ also define, through Theorem~\ref{thm:ACN}, a Sasakian Lie algebra with $\dd$ in place of $\h_1$.

\begin{lemma}\label{lem:data-d}
The K\"ahler Lie algebra $(\dd,[\cdot,\cdot]_\dd,J,\ip{\cdot}{\cdot})$, together with the same
$k$, $\gamma$ and $H_0$, satisfies the conditions \eqref{compat} of Theorem~\ref{thm:ACN}, with
the same one-form $\mu$. Consequently Theorem~\ref{thm:ACN} produces a Sasakian Lie algebra
$\g_0=\R\xi\oplus\dd\oplus\q$ with trivial centre and $\dim\Imm\ad_\xi=2$, which is solvable
because $H_0\neq0$ and $\dd$ is solvable.
\end{lemma}

\begin{proof}
By \cite[Proposition~6.4]{ACN}, a modification does not change the derived algebra, so
\[
\theta(\h_1,\h_1)=[\dd,\dd]=\n .
\]
In particular $H_0\in[\dd,\dd]$ and $JH_0=-a\,h_i\in\aaa=[\dd,\dd]^\perp$. Since the metric
and $J$ are those of $\h_1$, the one-form $\mu=-\tfrac12\ip{\cdot}{JH_0}$ is unchanged.

Next, for $x\in\dd$,
\[
\theta(x,H_0)=[x,H_0]_\dd+\phi_xH_0-\phi_{H_0}x=[x,H_0]_\dd ,
\]
because $\phi_xH_0=0$ by Proposition~\ref{prop:mod}(b), and $\phi_{H_0}=0$ since $\phi$
vanishes on $[\dd,\dd]\ni H_0$ (Definition~\ref{def:mod}). Hence $[x,H_0]_\dd=2\mu(x)H_0$.

Finally we check that $\omega=-\frac{1}{2k}\,d_\dd\gamma$, where $d_\dd$ denotes the
Chevalley--Eilenberg differential of $[\cdot,\cdot]_\dd$. By Proposition~\ref{prop:mod}(c),
$\omega=-d_\phi f=-d_\dd f$, while $\omega=-\frac1{2k}d_\phi\gamma$ by \eqref{compat}.
Therefore $d_\phi(\gamma-2kf)=0$, that is, $\gamma-2kf$ vanishes on
$\theta(\h_1,\h_1)=[\dd,\dd]$. But a one-form vanishing on $[\dd,\dd]$ is also
$d_\dd$-closed, so
\[
d_\dd\gamma=2k\,d_\dd f=-2k\,\omega ,
\]
as required.
\end{proof}

Let $\s_0\subset\g_0$ be the quotient model of $\g_0$ provided by Lemma~\ref{Lmap}, endowed
with the metric and complex structure of Lemma~\ref{Lisom}. Writing $L_0$ for the map of
Lemma~\ref{lem:L} attached to $\g_0$ and identifying $L_0(x)$ with $x\in\dd$, we have
$\s_0=\dd\oplus\q$ with brackets
\[
[x,y]_0=[x,y]_\dd, \quad [x,u]_0=\mu(x)u, \quad [x,v]_0=\mu(x)v, \quad [u,v]_0=H_0,
\qquad x,y\in\dd,
\]
the metric making $\dd\perp\q$ and $\{u,v\}$ orthonormal, and the complex structure being $J$
on $\dd$ and $u\mapsto v$ on $\q$. In particular $[\cdot,\cdot]_0$ is a Lie bracket, being
that of a subalgebra of $\g_0$, and Proposition~\ref{CHfactor} applies to $\s_0$. We identify $\s_0$ with $\s$ as vector spaces via $x\mapsto L(x)$
for $x\in\dd=\h_1$ and the identity on $\q$.

\begin{lemma}\label{lift}
For $x\in\dd$, let $\widehat\phi_x$ denote the extension of $\phi_x$ to $\s_0$ by zero on
$\q$. Then
\[
\Psi\colon\s_0\longrightarrow\Der_{u}(\s_0),\qquad
\Psi_{x+w}=\widehat\phi_x\qquad(x\in\dd,\ w\in\q),
\]
is a modification, and $\s=(\s_0)_\Psi$ as metric K\"ahler Lie algebras. Consequently,
\[
\nabla^\s R^\s=\nabla^{\s_0}R^{\s_0}
\]
under the common vector-space identification.
\end{lemma}

\begin{proof}

Let $X=x+w\in\s_0$ with $x\in\dd$ and $w\in\q$, so that $\Psi_X=\widehat\phi_x$. We first
check that $\widehat\phi_x\in\Der_{u}(\s_0)$.
\begin{itemize}
    \item Each $\widehat\phi_x$ is skew-symmetric and commutes with the complex structure, since
$\phi_x$ does and $\widehat\phi_x|_\q=0$ with $J^\s\q=\q$.
\item Each $\widehat\phi_x$ is a derivation of $[\cdot,\cdot]_0$. On $\dd$ this holds because
$\phi_x$ is a derivation of $\dd$. For the brackets involving $\q$, observe that, for $y\in\dd$,
\[
\mu(\phi_xy)=-\tfrac12\ip{\phi_xy}{JH_0}=\tfrac12\ip{y}{\phi_xJH_0}=0 ,
\]
because $\phi_xJH_0=J\phi_xH_0=0$. Hence, for $y\in\dd$ and $w\in\q$,
\[
\widehat\phi_x[y,w]_0=\widehat\phi_x\big(\mu(y)w\big)=0
=\mu(\phi_xy)w=[\phi_xy,w]_0+[y,\widehat\phi_xw]_0 ,
\]
and the derivation identity on $[u,v]_0=H_0$ reads $\phi_xH_0=0$, which holds.
\end{itemize}
Thus $\widehat\phi_x\in\Der_{\mathfrak u}(\s_0)$. We now check that $\Psi$ is a
modification:

\begin{itemize}
    \item $[\Psi_X,\Psi_Y]=0$ and $\Psi_{\Psi_XY}=0$ follow from the corresponding
conditions for $\phi$, since $\Psi_XY=\phi_xy\in\dd$.
\item $\Psi_{[X,Y]_0}=0$: writing
$X=x+w$, $Y=y+z$ with $x,y\in\dd$ and $w,z\in\q$, the $\dd$-component of
$[X,Y]_0$ is $[x,y]_\dd+tH_0$ for some $t\in\R$, so
\[
\Psi_{[X,Y]_0}=\widehat\phi_{[x,y]_\dd}+t\,\widehat\phi_{H_0}=0 ,
\]
using $\phi_{[x,y]_\dd}=0$ and $\phi_{H_0}=0$.
\end{itemize}
Hence $\Psi$ is a modification. Its modified bracket is
\[
[x,y]_\Psi=[x,y]_\dd+\phi_xy-\phi_yx=[x,y]_\phi\qquad \text{for } x,y\in\dd.
\]
The brackets involving $\q$ are unchanged, since $\widehat\phi$ vanishes on $\q$ and
$\Psi_w=0$ for $w\in\q$. These are exactly the brackets of $\s$, therefore
$\s=(\s_0)_\Psi$. Finally, by Proposition~\ref{prop:mod} and Corollary \ref{Dacts},
\[
\nabla^\s_XR^\s=\nabla^{\s_0}_XR^{\s_0}+\Psi_X\cdot R^{\s_0}=\nabla^{\s_0}_XR^{\s_0} .
\qedhere
\]

\end{proof}

\subsection{The main theorem}\label{proofmain}

We can now state and prove the main result of the section. 

\begin{theorem}\label{main}
Let $\g$ be a solvable Sasakian Lie algebra with trivial centre and assume $\dim\Imm\ad_\xi=2$.
Then
\[
\g \text{ is }\varphi\text{-symmetric}\iff
\nabla^{\h_1}R^{\h_1}=0 \ \text{ and }\ \nabla^{\h_1}\rho_0=0 .
\]
\end{theorem}

\begin{proof}
Since $\g$ is solvable, $H_0\neq0$ and $\h_1$ is solvable by \cite[Proposition 4.5]{ACN}.
It follows from \eqref{quotient} that
\begin{equation*}\label{gsymm}
\g \text{ is }\varphi\text{-symmetric}\iff \nabla^\s R^\s=0 .
\end{equation*}
Let us prove both implications:

\noindent($\Rightarrow$) If $\nabla^\s R^\s=0$, then Proposition \ref{curv} gives 
$\nabla^{\h_1}R^{\h_1}=0$ and $\nabla^{\h_1}\rho_0=0$.

\smallskip
\noindent($\Leftarrow$) Assume $\nabla^{\h_1}R^{\h_1}=0$ and $\nabla^{\h_1}\rho_0=0$, and
write $\h_1=\dd_\phi$. Then:
\begin{itemize}
\item $\dd$ is locally symmetric and $\nabla^{\dd}\rho_0=0$, by Lemma \ref{modbase};
\item hence $\n_{(\epsilon_i\pm\epsilon_r)/2}=\{0\}$ for every $r>i$, by Proposition
\ref{rootobstruction};
\item hence $\nabla^{\s_0}R^{\s_0}=0$, by Proposition \ref{CHfactor} applied to $\s_0$;
\item hence $\nabla^\s R^\s=0$, by Lemma \ref{lift}.
\end{itemize}
This shows that $\g$ is $\varphi$-symmetric.
\end{proof}

Theorem \ref{main} can be rephrased in terms of the normal $j$-algebra $\mathfrak d$ underlying $\mathfrak h_1$ and of the Kähler quotient. Recall that $H_0=-a\,x_i$ selects a distinguished root $\epsilon_i$, and that $\mathfrak d_i=\mathrm{span}\{h_i,x_i\}\oplus^\perp\mathfrak n_{\epsilon_i/2}$ and $\mathfrak k_i=\mathfrak d_i^\perp$.
 
\begin{theorem}\label{thm:geometric}
Let $\mathfrak g$ be a solvable Sasakian Lie algebra with trivial centre and assume
$\dim\operatorname{Im}\operatorname{ad}_\xi=2$, and write
$\mathfrak h_1=\mathfrak d_\phi$ with $(\mathfrak d,J,f)$ a normal
$j$-algebra and $H_0=-a\,x_i$. The following are equivalent:
\begin{enumerate}
\item[(i)] $\mathfrak g$ is $\varphi$-symmetric;
\item[(ii)] $\mathfrak d=\mathfrak d_i\oplus^\perp\mathfrak k_i$ is an orthogonal
sum of commuting ideals with $\mathfrak k_i$ locally symmetric.
\end{enumerate}
In this case $\mathfrak d_i$ is the solvable model of
$\mathbb{C}H^{n+1}(-a^2)$, where $2n=\dim\mathfrak n_{\epsilon_i/2}$, and the
Kähler quotient $G/H$ is locally isometric to the Hermitian symmetric space
\[
\mathbb{C}H^{n+2}(-a^2)\times K ,
\]
where $K$ is the Hermitian symmetric space of noncompact type whose solvable
model is $\mathfrak k_i$ (and $K$ is a point if $\mathfrak k_i=0$).
\end{theorem}
 
\begin{proof}
By Lemma \ref{modbase}, both the local symmetry of $\mathfrak h_1$ and the condition $\nabla^{\mathfrak h_1}\rho_0=0$ may be tested on $\mathfrak d$. It follows from 
Proposition \ref{rootobstruction} that $\nabla^{\mathfrak d}\rho_0=0$ holds if and only if
$\mathfrak n_{(\epsilon_i\pm\epsilon_r)/2}=\{0\}$ for every $r>i$. We now show that the vanishing of these mixed root spaces is equivalent to $\mathfrak d_i$ and $\mathfrak k_i$ being commuting ideals of $\mathfrak d$. One implication is Lemma \ref{lem:rho0-commuting-ideals}. For the converse, if $\mathfrak k_i$ is a Lie subalgebra and $0\neq x\in\mathfrak n_{(\epsilon_i-\epsilon_r)/2}\subseteq\mathfrak k_i$
for some $r>i$, then $Jx\in\mathfrak n_{(\epsilon_i+\epsilon_r)/2} \subseteq \mathfrak k_i$
and $\theta(x,Jx)\in\mathfrak n_{\epsilon_i}$, so by  \eqref{ftrick} we have 
$\theta(x,Jx)=a\,\omega(x,Jx)\,x_i=-a|x|^2x_i\notin\mathfrak k_i$, a
contradiction; hence $\mathfrak n_{(\epsilon_i-\epsilon_r)/2}=\{0\}$, and
$\mathfrak n_{(\epsilon_i+\epsilon_r)/2}=J\mathfrak n_{(\epsilon_i-\epsilon_r)/2}=\{0\}$.
 
Assume now that $\mathfrak d=\mathfrak d_i\oplus^\perp\mathfrak k_i$ is a sum
of commuting ideals. The brackets of $\mathfrak d_i$ are
\[
\theta(h_i,x_i)=a\,x_i,\qquad \theta(h_i,y)=\tfrac a2\,y,\qquad
\theta(x_i,y)=0,\qquad \theta(y,y')=a\,\omega(y,y')\,x_i,
\]
for any $y,y'\in\mathfrak n_{\epsilon_i/2}$, where the last equation follows  from \eqref{ftrick}. Hence, with $\hat H:=h_i$, $\hat Z:=-x_i$ and
$\mathfrak v:=\mathfrak n_{\epsilon_i/2}$, the relations \eqref{CHbrackets} hold with $\dim\mathfrak v=2n$, so $\mathfrak d_i$ is the solvable model of
$\mathbb{C}H^{n+1}(-a^2)$ by Remark \ref{rem:solvmodel}; in particular, $\mathfrak d_i$ is locally symmetric. Now,  $\mathfrak d$ is locally symmetric if and only if $\mathfrak k_i$ is. Putting everything together,
\[
\text{(i)}\iff\nabla^{\mathfrak h_1}R^{\mathfrak h_1}=0\ \text{and}\
\nabla^{\mathfrak h_1}\rho_0=0 \iff \text{(ii)},
\]
the first equivalence being Theorem \ref{main}.

For the last assertion,  it follows from Lemma \ref{lift} that the Kähler quotient $\mathfrak s$ of $\mathfrak g$ is a modification of the Kähler quotient $\mathfrak s_0$ of the Sasakian Lie algebra $\mathfrak g_0$ attached to $\mathfrak d$,
and we know that modifications preserve the curvature tensor (Proposition \ref{prop:mod}(a)); hence $\mathfrak s$ and $\mathfrak s_0$ have the same curvature and, being locally symmetric, their simply connected groups are locally isometric. By
Proposition \ref{CHfactor} applied to $\mathfrak s_0$,
$\mathfrak s_0=\mathfrak s_i\oplus^\perp\tilde{\mathfrak k}_i$ with $\mathfrak s_i$ the solvable model of $\mathbb{CH}^{n+2}(-a^2)$ and
$\tilde{\mathfrak k}_i$ isometric to $\mathfrak k_i$, whose simply connected
group is $K$ by Remark \ref{rem:solvmodel}. Since $p|_S\colon S\to G/H$ is a diffeomorphism and $S$ is the simply connected Lie group with Lie algebra $\mathfrak s$, the claim follows.
\end{proof}
 

\begin{remark}
Theorems \ref{H0zero} and \ref{main} leave out the case $H_0\neq0$ with $\mathfrak h_1$ not solvable, in which $\mathfrak g$ is not solvable by \cite[Proposition 4.5]{ACN}. We do not know whether nonsolvable Kähler Lie algebras exist. If they do not exist, Theorems \ref{H0zero} and \ref{main} cover every Sasakian Lie algebra with trivial centre and $\dim\operatorname{Im} \operatorname{ad}_\xi=2$.
\end{remark}
\section{Low dimensions and examples}\label{lowdim}

It follows from \cite{BV1} that the Sasakian Lie algebras $\mathfrak{su}(2),\, \mathfrak{sl}(2,\R)$ and $\h_3$ are $\varphi$-symmetric. We now study what happens in higher dimensions.

\begin{proposition}\label{dim7}
Every Sasakian Lie algebra $\g$ with trivial centre, $\dim\Imm\ad_\xi=2$ and $5\leq \dim\g\leq7$ is
$\varphi$-symmetric.
\end{proposition}

\begin{proof}
Here $\dim\h_1=\dim\g-3\in\{2,4\}$. If $\dim\h_1=2$ then $\h_1$ is the solvable model of
$\CH^1(-a^2)$ for some $a>0$. If $\dim\h_1=4$ then $\h_1$ is solvable by
\cite[Theorem 9]{Chu}, hence a modification of a normal $j$-algebra, which is
one of the two families of \cite[Theorem 5.9]{ACN}, namely $\dd_{a,b}$, the orthogonal
product of the solvable models of $\CH^1(-a^2)$ and $\CH^1(-b^2)$, or $\dd_a$, the solvable
model of $\CH^2(-a^2)$. In every case the underlying normal $j$-algebra is locally
symmetric, hence so is $\h_1$ by Lemma \ref{modbase}; moreover $\h_1$ is solvable.

Thus $\nabla^{\h_1}R^{\h_1}=0$ holds automatically in this dimension range, and only the
behaviour of $H_0$ remains to be examined. If $H_0=0$, local symmetry of $\h_1$ is all that
is required and Theorem \ref{H0zero} gives the claim. If $H_0\neq0$, then $\g$ is solvable
by \cite[Proposition 4.5]{ACN}, so Theorem \ref{main} applies and we are left with the
condition $\nabla^{\h_1}\rho_0=0$; by Proposition \ref{rootobstruction} this amounts to the
vanishing of the mixed root spaces, which is automatic here, since the underlying normal
$j$-algebra has either rank one, and then there are no mixed root spaces at all, or rank
two, and then $\aaa$, $\n_{\epsilon_1}$ and $\n_{\epsilon_2}$ already account for the four
dimensions of $\h_1$.
\end{proof}

\begin{remark}
For $\dim\g=5$, one has that $\dim\Imm\ad_\xi=2$ (the possibility $\dim \Imm\ad_\xi=4$ is excluded by \cite[Proposition 8.1]{ACN}); hence, Proposition \ref{dim7} recovers \cite[Theorem 5]{LL}.
\end{remark}

\begin{example}\label{easy}
Let $\dd$ be the normal $j$-algebra described by
\[
  \dd=\operatorname{span}\{h_1,h_2,\dots,h_n\}\oplus\operatorname{span}\{x_1,x_2, \dots,x_n\}
\]
with $\mathfrak a=\operatorname{span}\{h_1,h_2,\dots,h_n\}$ and root spaces
\[
  \mathfrak n_{\epsilon_i}=\mathbb R x_i,\qquad
  i=1,\dots,n,
\]
and there are no other roots; here $\epsilon_i=a_i\,h_i^{*}$,  with $a_i>0$, and $Jx_i=h_i$. The brackets are
\[
  \theta(h_i,x_i)=a_i\,x_i,\qquad i=1,\dots,n,
\]
and therefore 
\[
  \omega=\sum _i^n h_i^{*}\wedge x_i^{*}=-df \qquad \text{with} \qquad f=\sum _{i}^n\tfrac{1}{a_i}x_i^{*}.
\]
Thus $\dd$ is the orthogonal product of $n$ copies of $\mathfrak{aff}(\mathbb R)$, that
is, the solvable model of $\mathbb{CH}^1(-a_1^2)\times\mathbb{CH}^1(-a_2^2)\times \dots \times 
\mathbb{CH}^1(-a_n^2)$; in particular it is locally symmetric for all $a_i$.

Without loss of generality, we may take $H_0=-a_1\,x_1$. Setting also $k=1$ and $\gamma=2f$, we have that $JH_0=-a_1\,h_1$,
$\mu=\tfrac{a_1}{2}h_1^{*}$ and $\gamma(H_0)=-2$. There are no mixed root spaces at all,
so $\nabla^{\dd}\rho_0=0$ by Proposition \ref{rootobstruction}; concretely
\[
  \rho_0=-a_1^{2}\,h_1^{*}\wedge x_1^{*}=-a_1^{2}\,\omega|_{\dd_1}\oplus 0,
  \qquad \dd_1=\operatorname{span}\{h_1,x_1\},\quad \kk_1=\dd_1^{\perp}.
\]
By Theorem \ref{main} the $(2n+3)$-dimensional Sasakian Lie algebra
$\g=\mathbb R\xi\oplus\dd\oplus\q$ is $\varphi$-symmetric.
\end{example}

\medskip

If in the previous example we choose $n=3$, we obtain an example of a $9$-dimensional Sasakian Lie algebra that is $\varphi$-symmetric. In the next example we exhibit a $6$-dimensional symmetric normal $j$-algebra which gives rise to $9$-dimensional Sasakian Lie algebras that are not  $\varphi$-symmetric.

\begin{example}\label{falla-dim9}
Let $\dd$ be the normal $j$-algebra described by
\[
  \dd=\operatorname{span}\{h_1,h_2\}\oplus\operatorname{span}\{x_1,x_2, U,V\}
\]
with $\aaa=\operatorname{span}\{h_1,h_2\}$. The root spaces are
\[
  \mathfrak n_{\epsilon_1}=\mathbb R x_1,\qquad \mathfrak n_{\epsilon_2}=\mathbb R x_2,\qquad
  \mathfrak n_{(\epsilon_1-\epsilon_2)/2}=\mathbb R U,\qquad
  \mathfrak n_{(\epsilon_1+\epsilon_2)/2}=\mathbb R V .
\]
Here $\epsilon_j=ah_j^{*}$, with $Jx_j=h_j$ and $JU=V$, and
\[
  f=\tfrac1a\bigl(x_1^{*}+x_2^{*}\bigr),\qquad
  \omega=h_1^{*}\wedge x_1^{*}+h_2^{*}\wedge x_2^{*}-U^{*}\wedge V^{*}
\]
satisfies $\omega=-df$, and the brackets are
\[
\begin{aligned}
&\theta(h_1,x_1)=ax_1, &&\theta(h_2,x_2)=ax_2, &&\theta(U,V)=-ax_1,\\
&\theta(h_1,U)=\frac{a}{2}U, &&\theta(h_2,U)=-\frac{a}{2}U, &&\theta(h_1,V)= \frac{a}{2}V,\\
&\theta(h_2,V)=\frac{a}{2}V, &&\theta(U,x_2)=-aV. &&
\end{aligned}
\]
A direct computation from these Lie brackets shows that $\dd$ is locally symmetric.

Now set
\[
H_0=-ax_1,\qquad k=1,\qquad \gamma=2f,\qquad \mu=\frac{a}{2}h_1^* .
\]
The two-form is $\rho_0=-a^2\,h_1^*\wedge x_1^*+a^2\,U^*\wedge V^*$, and it is not parallel:
from $\rho_0(U,V)=a^2$, $\rho_0(h_2,\cdot)=0$ and
$\nabla^{\dd}_Uh_2=\nabla^{\dd}_{h_2}U+\theta(U,h_2)=\frac{a}{2}U$ by \eqref{nablaa}, we get
\[
\big(\nabla^{\dd}_U\rho_0\big)(h_2,V)
=-\frac{a}{2}\rho_0(U,V)-\rho_0\big(h_2,\nabla^{\dd}_UV\big)=-\frac{a^3}{2} .
\]
Hence
 $\g$ is not
$\varphi$-symmetric, although its K\"ahler base $\dd$ is symmetric.
\end{example}

\begin{remark}\label{dim9unique}
In dimension nine, the  K\"ahler base $(\h_1,H_0)$ of Example \ref{falla-dim9} is, up to scaling, the only six-dimensional Kähler base giving rise to a solvable Sasakian structure in our class whose base is locally symmetric but whose Sasakian structure is not $\varphi$-symmetric.
Indeed, $\dim\dd=6$, and such a failure requires a nonzero pair of mixed root spaces,
so that
\[
  6=\underbrace{2}_{\mathfrak a}
   +\underbrace{1+1}_{\mathfrak n_{\epsilon_1},\ \mathfrak n_{\epsilon_2}}
   +\underbrace{2}_{\mathfrak n_{(\epsilon_1-\epsilon_2)/2}\oplus \mathfrak n_{(\epsilon_1+\epsilon_2)/2}} .
\]
Hence $\mathfrak n_{\epsilon_j/2}=\{0\}$ for every $j$, and there is no third
distinguished root: the root data must be
those of Example \ref{falla-dim9}, and $H_0$ must sit at the first index by Lemma \ref{smaller}.
Writing $a=\epsilon_1(h_1)$ and $b=\epsilon_2(h_2)$, these data determine the brackets
\[
\begin{aligned}
&\theta(h_1,x_1)=a\,x_1, &&\theta(h_2,x_2)=b\,x_2, &&\theta(U,V)=-a\,x_1, &&\theta(U,x_2)=-b\,V,\\
&\theta(h_1,U)=\tfrac{a}{2}U, &&\theta(h_2,U)=-\tfrac{b}{2}U,
&&\theta(h_1,V)=\tfrac{a}{2}V, &&\theta(h_2,V)=\tfrac{b}{2}V,
\end{aligned}
\]
for which $\bigl(\nabla^{\dd}_U R^{\dd}\bigr)(h_2,V)U=\tfrac{3}{4}\,b\,(b^2-a^2)\,V$.
Hence $\dd$ is locally symmetric if and only if $a=b$, that is, if and only if it is
the algebra of Example \ref{falla-dim9} up to scaling. Finally, a direct computation
shows that this $\dd$ admits no nonzero skew-symmetric derivation commuting with $J$, so it
has no nontrivial modifications and $\h_1=\dd$.
\end{remark}

\medskip

In the next example we show that taking the product of the normal $j$-algebra in Example \ref{falla-dim9} with a 2-dimensional Kähler factor yields 11-dimensional $\varphi$-symmetric examples. 

\begin{example}\label{rank3}
Let $\dd$ now denote the following normal $j$-algebra:
\[
\dd=\spann\{h_1,h_2,h_3\}\oplus\spann\{x_1, x_2,x_3,  U, V\}
\]

\noindent with
$\aaa=\spann\{h_1,h_2,h_3\}$. The root spaces are
\[
\n_{\epsilon_1}=\R x_1,\qquad \n_{\epsilon_2}=\R x_2,\qquad \n_{\epsilon_3}=\R x_3,
\]
\[
\n_{(\epsilon_1-\epsilon_2)/2}=\R U,\qquad \n_{(\epsilon_1+\epsilon_2)/2}=\R V ,
\]
and there are no other roots. Set $a:=\epsilon_1(h_1)=\epsilon_2(h_2)>0$ and $b:=\epsilon_3(h_3)>0$, so that $\epsilon_1=a h_1^*$, $\epsilon_2=a h_2^*$ and $\epsilon_3=b h_3^*$, together with $Jx_j=h_j$ for $j=1,2,3$ and
$JU=V$. We have that
\[ \omega=h_1^*\wedge x_1^*+h_2^*\wedge x_2^*+h_3^*\wedge x_3^*-U^*\wedge V^*=-df, \qquad \text{with} \qquad 
f=\frac{1}{a}x_1^* +\frac{1}{a}x_2^* + \frac{1}{b}x_3^*.
\]
It can be shown that the nonzero brackets of $\dd$ are
\[
\begin{gathered}
\theta(h_1,x_1)=a x_1,\qquad \theta(h_2,x_2)=a x_2,\qquad \theta(h_3,x_3)=b x_3,\\
\theta(h_1,U)=\frac{a}{2}U,\qquad \theta(h_2,U)=-\frac{a}{2}U,\qquad
\theta(h_1,V)=\frac{a}{2}V,\qquad \theta(h_2,V)=\frac{a}{2}V,\\
\theta(U,V)=-a x_1,\qquad \theta(x_2,U)=a V.
\end{gathered}
\]

We note that $\dd$ is the orthogonal sum of the $J$-invariant commuting ideals
$\dd_1:=\spann\{h_1,x_1,h_2,x_2,U,V\}$ and $\dd_2:=\spann\{h_3,x_3\}$. It is clear that $\dd_2$ is locally symmetric (it is the solvable model of $\mathbb{C} H^1(-b^2)$); moreover, $\dd_1$ is the locally symmetric normal $j$-algebra from Example \ref{falla-dim9}. By the Iwasawa decomposition of noncompact symmetric spaces (see \cite{Hel}), this factor is the solvable model of the Hermitian symmetric space  $\operatorname{Sp}(2,\R)/\operatorname{U}(2)$, realized as the Siegel upper half-space.

We now construct the Sasakian Lie algebra $\g$ with $k=1$ and $\gamma=2f$, so that $\gamma(\theta(x,y))=2\omega(x,y)$
holds automatically from $\omega=-df$. For the distinguished  vector $H_0\neq 0$, it follows from Lemma \ref{smaller} that there are two possibilities: $H_0\in \n_{\epsilon_1}$ or $H_0\in \n_{\epsilon_3}$. We have:
\begin{itemize}
    \item[(a)] $H_0=-b x_3$: in this case the structure is $\varphi$-symmetric. Indeed,  $JH_0=-b h_3$ and $\mu=\frac{b}{2}h_3^*$. There are no mixed root spaces involving $\epsilon_3$, so
$\nabla^{\dd}\rho_0=0$ by Proposition \ref{rootobstruction}; explicitly,
$\rho_0=-b^2\,h_3^*\wedge x_3^*=-b^2\,\omega|_{\dd_2}\oplus0$. Since $\dd$ is locally symmetric, Theorem \ref{main} shows that the 11-dimensional Sasakian Lie algebra $\g=\R\xi \oplus \dd\oplus\q$ is $\varphi$-symmetric, and by Proposition \ref{CHfactor} its
K\"ahler quotient $\s$ is
\[
\CH^2(-b^2)\times \operatorname{Sp}(2,\R)/\operatorname{U}(2) .
\]

\item[(b)] $H_0=-ax_1$: in this case, the structure is not $\varphi$-symmetric. This follows from the same computations as in Example \ref{falla-dim9}.
\end{itemize}
\end{example}

\begin{remark}
Since in an irreducible Hermitian symmetric space of rank $r\ge2$ all the mixed root spaces are non-zero \cite[Ch. X, Table VI]{Hel}, Theorem \ref{thm:geometric} shows that such a base never carries a $\varphi$-symmetric structure with $H_0\neq0$; Example \ref{falla-dim9} is the smallest instance.
\end{remark}

\


\begin{thebibliography}{99}

\bibitem{ACN} 
A. Andrada, S. G. Chiossi and A. J. N\'u\~nez Sullca,
\emph{$\eta$-Einstein Sasakian Lie algebras}, Manuscripta Math. \textbf{177} (2026), Paper No.~12. 

\bibitem{AD}
A. Andrada and G. Dileo, \textit{Odd-dimensional counterparts of Abelian complex and hypercomplex structures}, Math. Nachr. \textbf{296} (2023), 470--508.

\bibitem{AFV} 
A. Andrada, A. Fino and L. Vezzoni,
\emph{A class of Sasakian $5$-manifolds}, Transform. Groups \textbf{14} (2009), 493--512.

\bibitem{Blair}
D.\ Blair, \textit{Riemannian geometry of contact and symplectic manifolds}, Progress in Math. 203, Birkhäuser (2010).

\bibitem{BV1}
D.\ Blair and L.\ Vanhecke, 
\textit{Symmetries and $\varphi$-symmetric spaces}, Tôhoku Math. J., II. Ser. \textbf{39} (1987), 373--383.

\bibitem{BV2}
D.\  Blair and L.\ Vanhecke, \textit{New characterization of $\varphi$-symmetric spaces}, Kodai Math. J. \textbf{10} (1987), 102--107.


\bibitem{BG}
C.\ Boyer and K. Galicki, \textit{Sasakian Geometry}, Oxford Mathematical Monographs, Oxford University Press (2008).

\bibitem{CF} 
G. Calvaruso and A. Fino, \emph{Five-dimensional $K$-contact Lie algebras}, Monatsh. Math. \textbf{167} (2012), 35--59.

\bibitem{Chu}
B.-Y. Chu, \textit{Symplectic homogeneous spaces}, Trans. Am. Math. Soc. \textbf{197} (1974), 145--159.

\bibitem{DAtri1980}
J.\ E.\ D'Atri, \emph{The curvature of homogeneous Siegel domains}, J. Differential Geom. \textbf{15} (1980), 61--70.

\bibitem{FI}
Th. Friedrich and S. Ivanov, \textit{Parallel spinors and connections with skew-symmetric torsion in string theory}, 
Asian J. Math. \textbf{6} (2002), 303--335.


\bibitem{He} E.\ Heintze, \emph{On homogeneous manifolds of negative curvature}, Math. Ann. \textbf{211} (1974), 23--34.

\bibitem{JK} J. A. Jim\'enez and O. Kowalski,
\emph{The classification of $\varphi$-symmetric Sasakian manifolds},
Monatsh. Math. \textbf{115} (1993), 83--98.

\bibitem{Hel}
S.\ Helgason, \textit{Differential Geometry, Lie Groups, and Symmetric Spaces}, Graduate Studies in Mathematics \textbf{34},  American Mathematical Society (2001).


\bibitem{KN2} 
S.\ Kobayashi and K.\ Nomizu, \emph{Foundations of Differential Geometry}, vol.~II, Interscience, New York, 1969.

\bibitem{KW} 
O.\ Kowalski and S.\ Wegrzynowski,
\emph{A classification of five-dimensional $\varphi$-symmetric spaces}, Tensor (N.S.) \textbf{46} (1987), 379--386.

\bibitem{LL} 
E. Loiudice and A. Lotta, \emph{On five dimensional Sasakian Lie algebras with trivial center}, Osaka J. Math. \textbf{55} (2018), 39--49.


\bibitem{Ok} 
M.\ Okumura, \emph{Some remarks on space with a certain contact structure}, T\^ohoku Math. J. (2) \textbf{14} (1962), 135--145.

\bibitem{PS}
I.\ Piatetski-Shapiro, \textit{
Automorphic functions and the geometry of classical domains}.
Mathematics and its Applications \textbf{8}. New York-London-Paris: Gordon and Breach Science Publishers  (1969).

\bibitem{Ta} 
T. Takahashi, \emph{Sasakian $\varphi$-symmetric spaces},
T\^ohoku Math. J. (2) \textbf{29} (1977), 91--113.


\end{thebibliography}
\end{document}